\documentclass[11pt, a4paper]{article}
\usepackage{latexsym,amssymb}
\usepackage{amsmath}
\usepackage[mathscr]{eucal}
\usepackage{color}
\usepackage[abbrev]{amsrefs}

\newtheorem{Theorem}{\bf Theorem}[section]
\newtheorem{Lemma}{\bf Lemma}[section]
\newtheorem{Proposition}{\bf Proposition}[section]
\newtheorem{Corollary}{\bf Corollary}[section]
\newtheorem{Remark}{\bf Remark}[section]
\newtheorem{Example}{\bf Example}[section]
\newtheorem{Definition}{\bf Definition}[section]

\newenvironment{theorem}{\begin{Theorem}$\!\!\!$}{\end{Theorem}}
\newenvironment{lemma}{\begin{Lemma}$\!\!\!$}{\end{Lemma}}
\newenvironment{proposition}{\begin{Proposition}$\!\!\!$}{\end{Proposition}}
\newenvironment{corollary}{\begin{Corollary}$\!\!\!$}{\end{Corollary}}
\newenvironment{remark}{\begin{Remark}$\!\!\!$}{\end{Remark}}

\newenvironment{definition}{\begin{Definition}$\!\!\!$}{\end{Definition}}

\numberwithin{equation}{section}

\numberwithin{equation}{section}

\begin{document}
\title{Existence of solutions to a fractional semilinear heat equation\\ 
in uniformly local weak Zygmund type spaces}
\author{Norisuke Ioku, Kazuhiro Ishige, and Tatsuki Kawakami}
\date{}
\maketitle
\begin{abstract}
In this paper we introduce uniformly local weak Zygmund type spaces, 
and obtain an optimal sufficient condition for the existence of solutions 
to the critical fractional semilinear heat equation.
\end{abstract}
\vspace{10pt}
\vspace{3pt}

\vspace{40pt}
\noindent 
Addresses:

\smallskip
\noindent 
N. I.: Mathematical Institute, Tohoku University,
Aoba, Sendai 980-8578, Japan\\
\noindent 
E-mail: {\tt ioku@tohoku.ac.jp}\\

\smallskip
\noindent 
K. I.: Graduate School of Mathematical Sciences, The University of Tokyo,\\ 
3-8-1 Komaba, Meguro-ku, Tokyo 153-8914, Japan\\
\noindent 
E-mail: {\tt ishige@ms.u-tokyo.ac.jp}\\

\smallskip
\noindent 
{T. K.}: Applied Mathematics and Informatics Course,\\ 
Faculty of Advanced Science and Technology, Ryukoku University,\\
1-5 Yokotani, Seta Oe-cho, Otsu, Shiga 520-2194, Japan\\
\noindent 
E-mail: {\tt kawakami@math.ryukoku.ac.jp}\\

\newpage
\section{Introduction}
Consider the Cauchy problem 
for the fractional semilinear heat equation 
\begin{equation}
\tag{P}
\label{eq:P}
\left\{
\begin{array}{ll}
\partial_t u+(-\Delta)^{\frac{\theta}{2}}u=|u|^{p-1}u, & x\in{\mathbb R}^n,\,\,t>0,\vspace{5pt}\\
u(x,0)=\varphi(x), & x\in{\mathbb R}^n,
\end{array}
\right.
\end{equation}
where $n\ge 1$, $\partial_t:=\partial/\partial t$, $\theta\in(0,2]$, $p>1$, and $\varphi$ is a locally integrable function in ${\mathbb R}^n$. 
Here $(-\Delta)^{\theta/2}$ denotes the fractional power of the Laplace operator $-\Delta$ in ${\mathbb R}^n$. 
In this paper we establish the local-in-time existence of solutions to problem~\eqref{eq:P} in the critical case 
$$
p=p_\theta:=1+\frac{\theta}{n}
$$
in the framework of uniformly local weak Zygmund type spaces. 

The solvability of the Cauchy problem for semilinear heat equations 
has fascinated many mathematicians since the pioneering work by Fujita~\cite{F}. 
The literature is very large and we just refer to the comprehensive monograph~\cite{QS} 
and the papers \cites{AD, BP, BC, FHIL01, FHIL02, FI01, FI02, GM, HI01, HIT02, IKK, IKO01, IKO02, KY, LS01, LS02, LRSV, M, RS, S, W, Z}, 
some of which are closely related to this paper and the others include recent developments in this subject. 
The study of the solvability of problem~\eqref{eq:P} is divided into the following three cases:
$$
\mbox{$1<p<p_\theta$ (subcritical case)};
\qquad
\mbox{$p>p_\theta$ (supercritical case)};
\qquad
\mbox{$p=p_\theta$ (critical case)}.
$$
We collect some known results on necessary conditions and sufficient conditions 
for the existence of solutions to problem~\eqref{eq:P}.
\begin{itemize}
\item[(1)] Subcritical case $(1<p<p_\theta)$
  \begin{itemize}
  \item[(a)] (Necessity)\\
  There exists $C_1=C_1(n,\theta,p)>0$ such that 
  if problem~\eqref{eq:P} possesses a nonnegative solution in ${\mathbb R}^n\times(0,T)$ for some $T>0$, 
  then
  $$
  \sup_{x\in{\mathbb R}^n}\int_{B(x,T^{1/\theta})}\varphi(y)\,dy \le C_1 T^{\frac{n}{\theta}-\frac{1}{p-1}}.
  $$
  See \cites{AD, BP} for $\theta=2$ and \cite{HI01} for $\theta\in(0,2]$. 
  \item[(b)] (Sufficiency)\\
  There exists $\epsilon_1=\epsilon_1(n,\theta,p)>0$ such that 
  if 
  $$
  \sup_{x\in{\mathbb R}^n}\int_{B(x,T^{1/\theta})}|\varphi(y)|\,dy\le \epsilon_1 T^{\frac{n}{\theta}-\frac{1}{p-1}}
  $$
  for some $T\in(0,\infty)$, then problem~\eqref{eq:P} possesses a solution in ${\mathbb R}^n\times(0,T)$. 
  (See e.g., \cites{AD, HI01, W}.)
  \end{itemize}
  The results in (a) and (b) (see also \eqref{eq:1.4}) 
  imply that, 
  for any nonnegative measurable initial function $\varphi$ in ${\mathbb R}^n$, 
  problem~\eqref{eq:P} possesses a local-in-time nonnegative solution if and only if 
  $$
  \sup_{x\in{\mathbb R}^n}\int_{B(x,1)}\varphi(y)\,dy<\infty.
  $$
\item[(2)] Supercritical case $(p>p_\theta)$
  \begin{itemize}
  \item[(a)] (Necessity)\\
  There exists $C_2=C_2(n,\theta,p)>0$ such that 
  if problem~\eqref{eq:P} possesses a nonnegative solution in ${\mathbb R}^n\times(0,T)$ for some $T>0$, 
  then
  $$
  \sup_{x\in{\mathbb R}^n}\sup_{\sigma\in(0,T^{1/\theta})}\,|B(x,\sigma)|^{\frac{\theta}{n(p-1)}-1}
  \int_{B(x,\sigma)}\varphi(y)\,dy \le C_2.
  $$
  See \cites{AD, BP} for $\theta=2$ and \cite{HI01} for $\theta\in(0,2]$. 
  \item[(b)] (Sufficiency)\\
  For any $r\in(1,\infty)$, 
  there exists $\epsilon_2=\epsilon_2(n,\theta, p,r)>0$ such that
  if 
  $$
  \sup_{x\in{\mathbb R}^n}\sup_{\sigma\in(0,T^{1/\theta})}\,|B(x,\sigma)|^{\frac{\theta}{n(p-1)}-\frac{1}{r}}
  \left[\int_{B(x,\sigma)}|\varphi(y)|^r\,dy\,\right]^{\frac{1}{r}}\le\epsilon_2
  $$
  for some $T\in(0,\infty]$, then problem~\eqref{eq:P} possesses a solution in ${\mathbb R}^n\times(0,T)$. 
  See \cites{KY, RS} for $\theta=2$ and \cites{HI01, IKO01, IKO02, Z} for $\theta\in(0,2]$. 
  (See e.g., \cites{AD, IKK, W} for related results.) 
  \end{itemize}
   \item[(3)] Critical case $(p=p_\theta)$
  \begin{itemize}
  \item[(a)] (Necessity)\\
  There exists $C_3=C_3(n,\theta)>0$ such that
  if problem~\eqref{eq:P} possesses a nonnegative solution in ${\mathbb R}^n\times(0,T)$ for some $T>0$, 
  then
  $$
  \sup_{x\in{\mathbb R}^n} \int_{B(x,\sigma)}\varphi(y)\,dy
  \le C_3\left[\log\biggr(e+\frac{T^{1/\theta}}{\sigma}\biggr)\right]^{-\frac{n}{\theta}},
  \quad \sigma\in(0,T^{1/\theta}).
  $$
  See \cites{BP} for $\theta=2$ and \cite{HI01} for $\theta\in(0,2]$.
  \item[(b)] (Sufficiency)\\
  For any $\alpha>0$, 
  there exists $\epsilon_3=\epsilon_3(n,\theta,\alpha)>0$ such that
  if 
  $$
  \sup_{x\in{\mathbb R}^n}\Psi_\alpha^{-1}\left[\,\frac{1}{|B(x,\sigma)|}\int_{B(x,\sigma)}
  \Psi_\alpha(T^\frac{1}{p-1}|\varphi(y)|)\,dy\,\right]\le\epsilon_3\rho(\sigma T^{-\frac{1}{\theta}}),
  \quad
  \sigma\in(0,T^{1/\theta}),
  $$
  for some $T>0$, 
  then problem~\eqref{eq:P} possesses a solution in ${\mathbb R}^n\times(0,T)$,
  where
  $$
  \Psi_\alpha(s):=s[\log (e+s)]^\alpha,
  \qquad
  \rho(s):=
  s^{-n}\biggr[\log\biggr(e+\frac{1}{s}\biggr)\biggr]^{-\frac{n}{\theta}}. 
  $$
  See \cites{HI01, IKO01, IKO02}.
  \end{itemize}
\end{itemize}
Furthermore, the results in (2) and (3) imply the following results.
\begin{itemize}
 \item[(4)] Let $p\ge p_\theta$, and set
  \begin{equation}
  \label{eq:1.1}
  \varphi_c(x):=
  \left\{
  \begin{array}{ll}
  |x|^{-n}\displaystyle{\biggr[\log\left(e+\frac{1}{|x|}\right)\biggr]^{-\frac{n}{\theta}-1}}
  \quad & \mbox{if}\quad \displaystyle{p=p_\theta},\vspace{7pt}\\
  |x|^{-\frac{\theta}{p-1}}\quad & \mbox{if}\quad \displaystyle{p>p_\theta},\vspace{3pt}
  \end{array}
  \right.
  \quad\mbox{for $x\in{\mathbb R}^n$}.
  \end{equation}
 \begin{itemize}
 \item[(a)] 
  There exists $C_4=C_4(n,\theta,p)>0$ such that if 
  $$
  \varphi(x)\ge C_4\varphi_c(x)
  $$ 
  for a.a.~$x$ in a neighborhood of the origin, then 
  problem~\eqref{eq:P} possesses no local-in-time nonnegative solutions.
  \item[(b)] 
  There exists $\epsilon_4=\epsilon_4(n,\theta,p)>0$ such that if 
  $$
  |\varphi(x)|\le\epsilon_4\varphi_c(x)+K,\quad\mbox{a.a.~$x\in{\mathbb R}^n$},  
  $$
  for some $K\ge 0$, then problem~\eqref{eq:P} possesses a local-in-time solution. 
  \end{itemize}
\end{itemize}
The results in (4) show that the ``strength" of the singularity at the origin of the function $\varphi_c$ 
is the critical threshold for the local-in-time solvability of problem (P). 
The function $\varphi_c$ is quite useful for identifying optimal function spaces to which initial functions belong 
from the view of the solvability of problem~\eqref{eq:P}.
We remark that assertion~(2)-(b) with $r=1$ and assertion~(3)-(b) with $\alpha=0$ do not hold.
(See \cite{T}*{Theorem~1, Proposition 1}, 
which treat only the case of $\theta=2$, but which is also applicable to the case of $\theta\in(0,2)$. See also \cite{KT}*{Section~4}.)

There are (at least) two useful strategies for the proof of the existence of solutions to problem~\eqref{eq:P}. 
One is the supersolution method (SSM) and the other is the contraction mapping theorem (CMT).
SSM depends on the following principle: if there exists a nonnegative supersolution $v$ 
to problem~\eqref{eq:P} in ${\mathbb R}^n\times(0,T)$ for some $T>0$, 
then problem~\eqref{eq:P} possesses a nonnegative solution $u$ in ${\mathbb R}^n\times(0,T)$ such that $u\le v$ in ${\mathbb R}^n\times(0,T)$. 
In our problem~\eqref{eq:P} with nonnegative initial function $\varphi$
the following functions have been used as supersolutions in ${\mathbb R}^n\times(0,T)$ for some $T>0$:
$$
2S_\theta(t)\varphi\quad (1<p<p_\theta);
\quad
2\left(S_\theta(t)\varphi^r\right)^{\frac{1}{r}}\quad (p>p_\theta);
\quad
2\Psi^{-1}_\alpha\left(S_\theta(t)\Psi_\alpha(\varphi)\right)\quad (p=p_\theta),
$$
where $S_\theta(t)\varphi$ is a solution to the fractional heat equation (see \eqref{eq:1.5}), 
$r>1$, and $\Psi_\alpha$ is as in assertion~(3)-(b). 
(See e.g., \cite{W} for $1<p<p_\theta$;  \cites{HI01,RS} for $p>p_\theta$; \cite{HI01} for $p=p_\theta$.)
Furthermore, thanks to the arguments in \cite{TW}, 
SSM is also applicable to the study of the existence of sign-changing solutions to problem~\eqref{eq:P} (see \cites{IKO01, IKO02}), 
however we require additional arguments in which sense the solution converges to the initial function. 
On the other hand, CMT is widely used in the proof of the existence of
solutions in various evolution equations, and the choice of function spaces is crucial. 
For our problem~\eqref{eq:P} with $p>p_\theta$, the existence of solutions has been proved by CMT 
in the framework of weak Lebesgue spaces (see \cites{FI01, IKK}) 
and Morrey spaces (see \cites{KY,Z}). 
The results in \cites{FI01, IKK, KY, Z} cover the result in (4)-(b) with $p>p_\theta$.
However, in the critical case $p=p_\theta$, 
the arguments in \cites{FI01, IKK, KY, Z} are not applicable to the proof of assertion~(4)-(b),
by the logarithmic singularity of $\varphi_c$.
$$
\begin{tabular}{|c||c||c|c|}
 \hline
  & Supersolution method (SSM) & Weak spaces (CMT) & Morrey spaces (CMT)\\
 \hline \hline
$p>p_\theta$ & See \cites{HI01, RS}  & See \cites{FI01, IKK} & See \cites{KY,Z} \\
 \hline
$p=p_\theta$
 & See \cite{HI01} & Open & Not applicable (see \cite{T}) \\
 \hline
\end{tabular}
$$

The aim of this paper is to establish a sharp sufficient condition
on the existence of solutions to problem~\eqref{eq:P} in the critical case $p=p_\theta$ in the framework of Banach spaces. 
For the critical case $p=p_\theta$, 
the weak Zygmund space $L^{1,\infty}(\log L)^{1+n/\theta}$ 
seems a reasonable Banach space since $\varphi_c\in L^{1,\infty}(\log L)^{1+n/\theta}$.
(See Remark~\ref{Remark:4.1}~(i) for the definition of the weak Zygmund spaces $L^{q,\infty}(\log L)^\alpha$, where $1\le q<\infty$ and $\alpha\ge 0$.)  
Then we require sharp decay estimates of solutions to the fractional heat equation in the weak Zygmund spaces $L^{q,\infty}(\log L)^\alpha$, 
however, by the peculiarity of $L^{1,\infty}(\log L)^\alpha$, 
it seems difficult to obtain our desired sharp decay estimates. 
(See Remark~\ref{Remark:4.1}~(ii) for further details.)

In this paper
we introduce new weak Zygmund type spaces ${\mathfrak L}^{q,\infty}(\log {\mathfrak L})^\alpha$ 
and uniformly local weak Zygmund type spaces ${\mathfrak L}_{{\rm ul}}^{q,\infty}(\log {\mathfrak L})^\alpha$. 
Then we establish sharp decay estimates of solutions to the fractional heat equation in the spaces ${\mathfrak L}^{q,\infty}(\log {\mathfrak L})^\alpha$ 
and ${\mathfrak L}_{{\rm ul}}^{q,\infty}(\log {\mathfrak L})^\alpha$, 
and obtain a sufficient condition on the existence of solutions to problem~\eqref{eq:P} with $p=p_\theta$ 
in the framework of the space ${\mathfrak L}_{{\rm ul}}^{q,\infty}(\log {\mathfrak L})^\alpha$. 
Our sufficient condition is simpler than that of assertion~(3)-(b) and covers assertion~(4)-(b) with $p=p_\theta$ . 
\vspace{5pt}

We introduce some notations and define the weak Zygmund type spaces ${\mathfrak L}^{q,\infty}(\log {\mathfrak L})^\alpha$ 
and the uniformly local weak Zygmund type spaces ${\mathfrak L}_{{\rm ul}}^{q,\infty}(\log {\mathfrak L})^\alpha$.
We also formulate the definition of solutions to problem~\eqref{eq:P}. 
Let ${\mathcal M}$ be the set of Lebesgue measurable sets in ${\mathbb R}^n$. 
For any $E\in{\mathcal M}$, we denote by $|E|$ and $\chi_E$ the $n$-dimensional Lebesgue measure of $E$ 
and the characteristic function of $E$, respectively. 
Let $L^1_{{\rm loc}}$ be the set of locally integrable functions in ${\mathbb R}^n$.
For any $q\in[1,\infty]$, we denote by $L^q$ and $\|\cdot\|_{L^q}$ the usual $L^q$-space on ${\mathbb R}^n$ and its norm, respectively.

Let $q\in[1,\infty]$ and $\alpha\in[0,\infty)$.
We define the weak Zygmund type space ${\mathfrak L}^{q,\infty}(\log {\mathfrak L})^\alpha$ by
$$
	{\mathfrak L}^{q,\infty}(\log {\mathfrak L})^\alpha
	:=\{f\in L^1_{{\rm loc}}\,:\,\|f\|_{{\mathfrak L}^{q,\infty}(\log{\mathfrak L})^\alpha}<\infty\},
$$
where 
\begin{equation}
\label{eq:1.2}
\|f\|_{{\mathfrak L}^{q,\infty}(\log{\mathfrak L})^\alpha}
:=
\left\{
\begin{array}{ll}
\displaystyle{\sup_{s>0}\,\left\{\left[\log\left(e+\frac{1}{s}\right)\right]^\alpha
\sup_{|E|=s}\int_E |f(x)|^q\,dx\right\}^{\frac{1}{q}}} 
& \mbox{if}\,\,\, q\in[1,\infty),
\vspace{7pt}
\\
\|f\|_{L^\infty} 
& \mbox{if}\,\,\, q=\infty.\vspace{3pt}
\end{array}
\right.
\end{equation}
Then ${\mathfrak L}^{q,\infty}(\log{\mathfrak L})^\alpha$ is a Banach space equipped with the norm~$\|\cdot\|_{{\mathfrak L}^{q,\infty}(\log{\mathfrak L})^\alpha}$
(see Lemma~\ref{Lemma:2.1}). 
See \eqref{eq:2.9} for different expressions of the norm $\|\cdot\|_{{\mathfrak L}^{q,\infty}(\log{\mathfrak L})^\alpha}$. 
Remark that
\begin{equation}
\label{eq:1.3}
L^q
=
{\mathfrak L}^{q,\infty}(\log{\mathfrak L})^0
\supset
{\mathfrak L}^{q,\infty}(\log{\mathfrak L})^{\alpha}
\quad\mbox{for}\quad \alpha\ge 0. 
\end{equation}

Next, we introduce the uniformly local weak Zygmund type space ${\mathfrak L}_{{\rm ul}}^{q,\infty}(\log{\mathfrak L})^\alpha$ by
$$
{\mathfrak L}_{{\rm ul}}^{q,\infty}(\log{\mathfrak L})^\alpha:=\{f\in L^1_{{\rm loc}}\,:\,\|f\|_{\mathfrak L_{{\rm {ul}}}^{q,\infty}(\log{\mathfrak L})^\alpha}<\infty\},
$$
where 
$$
\|f\|_{{\mathfrak L}_{{\rm ul}}^{q,\infty}(\log{\mathfrak L})^\alpha}:=\sup_{z\in{\mathbb R}^n}\,\|f\chi_{B(z,1)}\|_{{\mathfrak L}^{q,\infty}(\log{\mathfrak L})^\alpha}.
$$
Then ${\mathfrak L}_{{\rm ul}}^{q,\infty}(\log{\mathfrak L})^\alpha$ is also a Banach space equipped with the norm~$\|\cdot\|_{{\mathfrak L}_{{\rm ul}}^{q,\infty}(\log{\mathfrak L})^\alpha}$.
For any $f\in {\mathfrak L}_{{\rm ul}}^{q,\infty}(\log{\mathfrak L})^\alpha$ and $\rho>0$, we often write
$$
|||f|||_{q,\alpha;\rho}:=\sup_{z\in{\mathbb R}^n}\,\|f\chi_{B(z,\rho)}\|_{{\mathfrak L}^{q,\infty}(\log{\mathfrak L})^\alpha}
$$
for simplicity. 
We remark that ${\mathfrak L}_{{\rm ul}}^{\infty,\infty}(\log{\mathfrak L})^\alpha=L^\infty$ and 
$|||\cdot|||_{\infty,\alpha;\rho}=\|\cdot\|_{L^\infty}$ for all $\alpha\in[0,\infty)$.
Notice that, for any $k\ge 1$, there exists $C=C(n,k)>0$ such that 
\begin{equation}
\label{eq:1.4}
|||f|||_{q,\alpha;k\rho}\le C|||f|||_{q,\alpha;\rho}
\end{equation}
for $f\in {\mathfrak L}_{{\rm ul}}^{q,\infty}(\log{\mathfrak L})^\alpha$ and $\rho>0$.

We formulate the definition of solutions to problem~\eqref{eq:P}. 
Let $\theta\in(0,2]$. Let $G_\theta$ be the fundamental solution to the fractional heat equation
$$
\partial_t v+(-\Delta)^{\frac{\theta}{2}}v=0\quad\mbox{in}\quad{\mathbb R}^n\times(0,\infty). 
$$
For any $\varphi$ in $L^1_{{\rm loc}}$,
we write
\begin{equation}
\label{eq:1.5}
(S_\theta(t)\varphi)(x):=\int_{{\mathbb R}^n}G_\theta(x-y,t)\varphi(y)\,dy,
\quad (x,t)\in{\mathbb R}^n\times(0,\infty),
\end{equation}
for simplicity.

\begin{definition}
\label{Definition:1.1}
Let $\theta\in(0,2]$, $p>1$, and $T>0$. Set $F_p(s):=|s|^{p-1}s$ for $s\in{\mathbb R}$. 
Let $u$ be a measurable and finite almost everywhere function in ${\mathbb R}^n\times(0,T)$. 
We say that $u$ is a solution to problem~\eqref{eq:P} in ${\mathbb R}^n\times(0,T)$ if, 
for a.a.~$(x,t)\in{\mathbb R}^n\times(0,T)$, 
\begin{itemize}
  \item 
  $G_\theta(x-y,t)\varphi(y)$ is integrable in ${\mathbb R}^n$ with respect to $y\in{\mathbb R}^n$;
  \item 
  $G_\theta(x-y,t-s)F_p(u(y,s))$ is integrable in ${\mathbb R}^n\times(0,t)$ with respect to $(y,s)\in{\mathbb R}^n\times(0,t)$;
  \item 
  $u$ satisfies
  $$
  u(x,t)=
  \left[S_\theta(t)\varphi\right](x)+\int_0^t \left[S_\theta(t-s)F_p(u(s))\right](x)\,ds.
  $$
\end{itemize}
\end{definition}

We are ready to state our main results. 
\begin{theorem}
\label{Theorem:1.1}
Let $\theta\in(0,2]$, $p=p_\theta= 1+\theta/n$, and $T_*\in(0,\infty)$. 
Then there exists $\epsilon>0$ such that if $\varphi\in {\mathfrak L}_{{\rm ul}}^{1,\infty}(\log{\mathfrak L})^{n/\theta}$ satisfies
$$
|||\varphi|||_{1,n/\theta;T^{1/\theta}}\le\epsilon\quad\mbox{for some $T\in(0,T_*]$},
$$
then problem~\eqref{eq:P} possesses a solution 
$u\in C((0,T):\mathfrak L^{1,\infty}_{{\rm ul}}(\log{\mathfrak L})^{n/\theta})
\cap L^\infty_{\rm loc}((0,T):L^\infty)$ in ${\mathbb R}^n\times(0,T)$, 
with~$u$ satisfying
\begin{equation}
\label{eq:1.6}
	\sup_{t\in(0,T)}\,|||u(t)|||_{1,n/\theta;T^{1/\theta}}
	+\sup_{t\in(0,T)}\,t^{\frac{n}{\theta}}\bigg[\log\bigg(e+\frac{1}{t}\bigg)\bigg]^{\frac{n}{\theta}}
	\|u(t)\|_{L^\infty}<\infty
\end{equation}
for $t\in(0,T)$. 
Furthermore, the solution~$u$ satisfies 
\begin{equation}
\label{eq:1.7}
	\begin{aligned}
	&
	\lim_{t\to +0}\|u(t)-S_\theta(t)\varphi\|_{{\mathfrak L}^{1,\infty}_{\rm ul}(\log{\mathfrak L})^\gamma}=0
	\qquad\mbox{for any $\gamma\in[0,n/\theta)$},
	\\
	&
	\lim_{t\to +0}u(t)=\varphi\quad\mbox{in the sense of distributions}.
	\end{aligned}
\end{equation}
\end{theorem}
We remark that Theorem~\ref{Theorem:1.1} with $T_*=\infty$ does not hold since problem~\eqref{eq:P} 
possesses no global-in-time positive solutions (see \cite{S}). 
As a direct consequence of Theorem~\ref{Theorem:1.1}, we obtain assertion~(4)-(b).
\begin{corollary}
\label{Corollary:1.1}
Let $\theta\in(0,2]$ and $p=p_\theta$. 
Let $\varphi_c$ be as in \eqref{eq:1.1}.
Then there exists $\epsilon>0$ such that if 
$$
|\varphi(x)|\le\epsilon\varphi_c(x)+K, \quad\mbox{a.a.~$x\in{\mathbb R}^n$},
$$
for some $K\ge 0$, then problem~\eqref{eq:P} possesses a local-in-time solution.
\end{corollary}
Furthermore, as a consequence of Theorem~\ref{Theorem:1.1}, we have:
\begin{theorem}
\label{Theorem:1.2}
Let $\theta\in(0,2]$ and $p=p_\theta$.
If $\varphi\in {\mathfrak L}_{{\rm ul}}^{1,\infty}(\log{\mathfrak L})^\alpha$ for some $\alpha>n/\theta$, then 
problem~\eqref{eq:P} possesses a solution~$u$ in ${\mathbb R}^n\times(0,T)$ for some $T>0$, 
with $u$ satisfying \eqref{eq:1.6} and \eqref{eq:1.7}.
\end{theorem}

The rest of this paper is organized as follows. 
In Section~2 we collect some properties of non-increasing rearrangement of measurable functions,  
and prove some lemmas in ${\mathfrak L}^{q,\infty}(\log{\mathfrak L})^\alpha$ and ${\mathfrak L}_{{\rm ul}}^{q,\infty}(\log{\mathfrak L})^\alpha$.
Furthermore, we recall Hardy's inequalities and some properties of $S_\theta(t)\varphi$. 
In Section~3  we establish decay estimates of $S_\theta(t)\varphi$ in weak Zygmund type spaces (see Proposition~\ref{Proposition:3.1}). 
Furthermore, we obtain decay estimates of $S_\theta(t)\varphi$ in ${\mathfrak L}_{{\rm ul}}^{q,\infty}(\log{\mathfrak L})^\alpha$ 
using Besicovitch's covering lemma. 
In Section~4 we apply the contraction mapping theorem in ${\mathfrak L}_{{\rm ul}}^{q,\infty}(\log{\mathfrak L})^\alpha$
to prove Theorems~\ref{Theorem:1.1} and \ref{Theorem:1.2}. 
In Appendix we give two propositions on relations 
among the weak Zygmund type spaces
$\mathfrak L^{q,\infty}(\log \mathfrak L)^\alpha$, the Zygmund spaces $L^q(\log L)^\alpha$, 
and the weak Zygmund spaces $L^{q,\infty}(\log L)^\alpha$. 
\section{Preliminaries}
In this section we introduce some notations, and give some lemmas on our weak Zygmund type spaces. 
Furthermore, we recall some lemmas on Hardy's inequalities. 
In all that follows we will use $C$ to denote generic positive constants and point out that $C$  
may take different values  within a calculation.
For any positive functions $f_1$ and $f_2$ in $(0,\infty)$, we write 
$$
\mbox{$f_1\asymp f_2$ for $s>0$}\quad\mbox{if}\quad
\mbox{$C^{-1}f_2(s)\le f_1(s)\le Cf_2(s)$ for $s>0$}.
$$
\subsection{Weak Zygmund type spaces}
For any (Lebesgue) measurable function $f$ in ${\mathbb R}^n$, 
we denote by $\mu_f$ the distribution function of~$f$, that is, 
$$
\mu_f(\lambda):=\left|\{x\,:\,|\,f(x)|>\lambda\}\right| 
\quad \mbox{for}\quad \lambda > 0.
$$
We define the non-increasing rearrangement $f^*$ of $f$ by 
$$
f^{*}(s):=\inf\{\lambda>0\,:\,\mu_f(\lambda)\le s\}\quad\mbox{for}\quad s\in[0,\infty).
$$
Here we adopt the convention $\inf\emptyset=\infty$. 
Then $f^*$ is non-increasing and right continuous in $[0,\infty)$, and it has the following properties (see \cite{Grafakos}*{Proposition~1.4.5}):
\begin{equation}
\label{eq:2.1}
(kf)^*=|k|f^*,\quad (|f|^q)^*=(f^*)^q,\quad 
\int_{{\mathbb R}^n}|f(x)|^q\,dx=\int_0^\infty f^*(s)^q\,ds,\quad
f^*(0)=\|f\|_{L^\infty},
\end{equation}
where $q\in(0,\infty)$ and $k\in{\mathbb R}$. 
We remark that if $E\in{\mathcal M}$ with $|E|<\infty$, then 
\begin{equation}
\label{eq:2.2}
(\chi_E)^*(s)=\chi_{[0,|E|)}(s)\quad\mbox{for}\quad s\ge 0.
\end{equation}
Define 
\begin{equation}
\label{eq:2.3}
f^{**}(s):=\frac{1}{s}\int_0^s f^*(\tau)\,d\tau\quad\mbox{for}\quad s\in(0,\infty).
\end{equation}
Here we collect properties of $f^*$ and $f^{**}$ used in the paper.
\begin{itemize}
  \item[(a)] 
  Since $f^*$ is non-increasing in $(0,\infty)$, it follows that 
  \begin{equation}
  \label{eq:2.4}
  f^{**}(s)\ge f^*(s)\quad\mbox{for}\quad s\in(0,\infty).
  \end{equation}
  \item[(b)] 
  For any $q\in[1,\infty)$, 
  Jensen's inequality together with \eqref{eq:2.1} yields
  \begin{equation}
  \label{eq:2.5}
  (f^{**}(s))^q\le \frac{1}{s}\int_0^s f^*(s)^q\,ds=\frac{1}{s}\int_0^s (|f|^q)^*(s)=(|f|^q)^{**}(s)\quad\mbox{for}\quad s\in(0,\infty).
  \end{equation} 
  \item[(c)] 
  It follows from \cite{BS}*{Chapter 2, Proposition 3.3} that 
  \begin{equation}
  \label{eq:2.6}
  f^{**}(s)=\frac{1}{s}\int_0^s f^*(\tau)\,d\tau=\frac{1}{s}\sup_{|E|=s}\int_E |f(x)|\,dx\quad\mbox{for}\quad s\in(0,\infty).
  \end{equation}
  \item[(d)]
  {(O'Neil's inequality)} For any $f$, $g\in L^1$, it follows from \cite{ONeil}*{Lemma~1.6} that 
  \begin{equation}
  \label{eq:2.7}
  (f*g)^{**}(s)\le\int_s^\infty  f^{**}(\tau)g^{**}(\tau)\,d\tau\quad\mbox{for}\quad s\in(0,\infty),
  \end{equation}
  where
  $$
  (f*g)(x)=\int_{{\mathbb R}^n}f(x-y)g(y)\,dy.
  $$
  \item[(e)]
  For any $f_1$, $f_2\in L^1_{{\rm loc}}$,
  it follows from \cite{ONeil}*{Theorem 3.3} that
  \begin{equation}
  \label{eq:2.8}
  (f_1f_2)^{**}(s)\le\frac{1}{s}\int_0^s f_1^*(\tau)f_2^*(\tau)\,d\tau\quad\mbox{for}\quad s\in(0,\infty).
  \end{equation}
\end{itemize}
Let $q\in[1,\infty)$ and $\alpha\ge 0$.
For any $L^1_{{\rm loc}}$-function $f$, 
by \eqref{eq:1.2}, \eqref{eq:2.1}, and \eqref{eq:2.6}
we have
\begin{equation}
\label{eq:2.9}
\begin{split}
\|f\|_{{\mathfrak L}^{q,\infty}(\log{\mathfrak L})^\alpha} & =\sup_{s>0}\,\left\{\left[\log\left(e+\frac{1}{s}\right)\right]^\alpha s(|f|^q)^{**}(s)\right\}^{\frac{1}{q}}\\
 & =\sup_{s>0}\,\left\{\left[\log\left(e+\frac{1}{s}\right)\right]^\alpha \int_0^s (|f|^q)^*(\tau)\,d\tau\right\}^{\frac{1}{q}}\\
 & =\sup_{s>0}\,\left\{\left[\log\left(e+\frac{1}{s}\right)\right]^\alpha \int_0^s f^*(\tau)^q\,d\tau\right\}^{\frac{1}{q}}.
\end{split}
\end{equation}
Furthermore, 
for any $E\in{\mathcal M}$ with $|E|<\infty$, 
it follows from \eqref{eq:2.2} and \eqref{eq:2.8} that 
\begin{equation}
\label{eq:2.10}
	(f\chi_E)^{**}(s)=(f\chi_E\chi_E)^{**}(s)
	\le\frac{1}{s}\int_0^s (f\chi_E)^*(\tau)(\chi_E)^*(\tau)\,d\tau
	=\frac{1}{s}\int_0^{\min\{s,|E|\}} (f\chi_E)^*(\tau)\,d\tau.
\end{equation}
For any $\beta\in[\alpha,\infty)$, since
\begin{equation}
\label{eq:2.11}
\mbox{the function $\displaystyle{(0,\infty)\ni \tau\mapsto \left[\log\left(e+\frac{1}{\tau}\right)\right]^{\alpha-\beta}\in{\mathbb R}}$ is non-decreasing},
\end{equation}
by \eqref{eq:2.1}, \eqref{eq:2.9}, and \eqref{eq:2.10} we have
\begin{align*}
	\|f\chi_E\|_{{\mathfrak L}^{1,\infty}(\log{\mathfrak L})^\alpha}
	& 
	=\sup_{s>0}\,\left\{\left[\log\left(e+\frac{1}{s}\right)\right]^\alpha s (f\chi_E)^{**}(s)\right\}\\
	&
	\le \sup_{s>0}\,\left\{\left[\log\left(e+\frac{1}{s}\right)\right]^\alpha \int_0^{\min\{s,|E|\}} (f\chi_E)^*(\tau)\,d\tau\right\}\\
 	& 
 	=\sup_{0<s\le |E|}\left\{\left[\log\left(e+\frac{1}{s}\right)\right]^{\beta+\alpha-\beta}
	\int_0^{\min\{s,|E|\}} (f\chi_E)^*(\tau)\,d\tau\right\}\\
 	& 
 	\le\left[\log\left(e+\frac{1}{|E|}\right)\right]^{\alpha-\beta}
 	\sup_{0<s\le |E|}\left\{\left[\log\left(e+\frac{1}{s}\right)\right]^{\beta}
	\int_0^{\min\{s,|E|\}} (f\chi_E)^*(\tau)\,d\tau\right\}\\
	& 
 	\le\left[\log\left(e+\frac{1}{|E|}\right)\right]^{\alpha-\beta}
 	\sup_{s>0}\left\{\left[\log\left(e+\frac{1}{s}\right)\right]^{\beta}
	\int_0^s (f\chi_E)^*(\tau)\,d\tau\right\}\\
 	& 
 	=\left[\log\left(e+\frac{1}{|E|}\right)\right]^{\alpha-\beta}
	\|f\chi_E\|_{{\mathfrak L}^{1,\infty}(\log{\mathfrak L})^\beta}.
\end{align*}
In particular, 
\begin{equation}
\label{eq:2.12}
|||f|||_{1,\alpha;\rho}\le C\left[\log\left(e+\frac{1}{\rho}\right)\right]^{\alpha-\beta}|||f|||_{1,\beta;\rho}
\end{equation}
for $f\in {\mathfrak L}_{{\rm ul}}^{1,\infty}(\log{\mathfrak L})^\beta$, $0\le\alpha\le\beta$, and $\rho>0$. 
Here we show that $\mathfrak L^{q,\infty}(\log \mathfrak L)^\alpha$ and $\mathfrak L_{\rm ul}^{q,\infty}(\log \mathfrak L)^\alpha$ are Banach spaces. 
\begin{lemma}
\label{Lemma:2.1}
For any $1\le q <\infty$ and $\alpha\ge 0$,
the weak Zygmund type space
$\mathfrak L^{q,\infty}(\log \mathfrak L)^\alpha$
and the uniformly local weak Zygmund type space
$\mathfrak L_{\rm ul}^{q,\infty}(\log \mathfrak L)^\alpha$ are Banach spaces.
\end{lemma}
{\bf Proof.}
Let $1\le q <\infty$ and $\alpha\ge 0$. 
It suffices to prove that $\mathfrak L^{q,\infty}(\log \mathfrak L)^\alpha$ 
(resp.~$\mathfrak L_{\rm ul}^{q,\infty}(\log \mathfrak L)^\alpha$) is a complete metric space
with the norm $\|\cdot\|_{\mathfrak L^{q,\infty}(\log \mathfrak L)^\alpha}$ (resp.~$\|\cdot\|_{\mathfrak L^{q,\infty}_{{\rm ul}}(\log \mathfrak L)^\alpha}$).  
Let $\{f_n\}$ be a Cauchy sequence in $\mathfrak L^{q,\infty}(\log \mathfrak L)^\alpha$.
It follows from \eqref{eq:1.3} that
$\{f_n\}$ is a Cauchy sequence in $L^{q}$, and
hence there exists $f\in L^q$ such that $f_n\to f$ as $n\to \infty$ in $L^q$.
Since the Cauchy sequence $\{f_n\}$ is bounded in 
$\mathfrak L^{q,\infty}(\log \mathfrak L)^\alpha$,
we observe from \eqref{eq:1.2} that $f\in \mathfrak L^{q,\infty}(\log \mathfrak L)^\alpha$. 
It remains to prove that $f_n\to f$ as $n\to \infty$
in $\mathfrak L^{q,\infty}(\log \mathfrak L)^\alpha$. 
For this aim, we take 
a subsequence 
$\{f_{n_k}\}$ such that $\{f_{n_k}\}$ converges almost everywhere to $f$.
Then Fatou's lemma gives us that
$$
\begin{aligned}
\|f-f_{n_j}\|_{\mathfrak L^{q,\infty}(\log \mathfrak L)^\alpha}
&
\le
\sup_{s>0}\left\{\left[\log\left(e+\frac{1}{s}\right)\right]^{\alpha}
\sup_{|E|=s}
\liminf_{k\to \infty}
\int_E |f_{n_k}(x)-f_{n_j}(x)|^q dx\right\}^{\frac{1}{q}}
\\
&
\le
\liminf_{k\to \infty}
\|f_{n_k}-f_{n_j}\|_{\mathfrak L^{q,\infty}(\log \mathfrak L)^\alpha}.
\end{aligned}
$$
This implies that $f_{n_j}$ converges to $f$ in $\mathfrak L^{q,\infty}(\log \mathfrak L)^\alpha$. 
Thus $\mathfrak L^{q,\infty}(\log \mathfrak L)^\alpha$ is a complete metric space.
Similarly, we see that $\mathfrak L_{\rm ul}^{q,\infty}(\log \mathfrak L)^\alpha$ is a complete metric space. 
The proof is complete.
$\Box$\vspace{5pt}

Next, we prove two lemmas on our weak Zygmund type spaces and uniformly local weak Zygmund type spaces. 
\begin{lemma}
\label{Lemma:2.2}
Let $q_1$, $q_2\in[1,\infty]$ and $\alpha_1$, $\alpha_2\ge 0$ be such that 
\begin{equation}
\label{eq:2.13}
1=\frac{1}{q_1}+\frac{1}{q_2},\quad \alpha=\frac{\alpha_1}{q_1}+\frac{\alpha_2}{q_2}.
\end{equation}
Then 
\begin{equation}
\label{eq:2.14}
\|f_1f_2\|_{{\mathfrak L}^{1,\infty}(\log{\mathfrak L})^\alpha}\le\|f_1\|_{{\mathfrak L}^{q_1,\infty}(\log{\mathfrak L})^{\alpha_1}}\|f_2\|_{{\mathfrak L}^{q_2,\infty}(\log{\mathfrak L})^{\alpha_2}}
\end{equation}
for $f_1\in {\mathfrak L}^{q_1,\infty}(\log{\mathfrak L})^{\alpha_1}$ and $f_2\in {\mathfrak L}^{q_2,\infty}(\log{\mathfrak L})^{\alpha_2}$. 
Furthermore, 
\begin{equation}
\label{eq:2.15}
|||\tilde{f}_1\tilde{f}_2|||_{1,\alpha;\rho}\le |||\tilde{f}_1|||_{q_1,\alpha_1;\rho} |||\tilde{f}_2|||_{q_2,\alpha_2;\rho}
\end{equation}
for $\tilde{f}_1\in {\mathfrak L}_{{\rm ul}}^{q_1,\infty}(\log{\mathfrak L})^{\alpha_1}$, $\tilde{f}_2\in {\mathfrak L}_{{\rm ul}}^{q_2,\infty}(\log{\mathfrak L})^{\alpha_2}$, and $\rho>0$.
\end{lemma}
{\bf Proof.}
Let $q_1$, $q_2\in[1,\infty)$ and $\alpha_1$, $\alpha_2\ge 0$ satisfy \eqref{eq:2.13}. 
Let $f_1\in {\mathfrak L}^{q_1,\infty}(\log{\mathfrak L})^{\alpha_1}$ and $f_2\in {\mathfrak L}^{q_2,\infty}(\log{\mathfrak L})^{\alpha_2}$. 
It follows from H\"older's inequality, \eqref{eq:2.8}, and \eqref{eq:2.9} that
\begin{align*}
 & \|f_1f_2\|_{{\mathfrak L}^{1,\infty}(\log{\mathfrak L})^\alpha}\\
 & =\sup_{s>0}\,\left\{\left[\log\left(e+\frac{1}{s}\right)\right]^\alpha s (f_1f_2)^{**}(s)\right\}\\
 & \le\sup_{s>0}\,\left\{\left[\log\left(e+\frac{1}{s}\right)\right]^\alpha\int_0^s f_1^*(\tau) f_2^*(\tau)\,d\tau\right\}\\
 & \le\sup_{s>0}\,\left\{\left[\log\left(e+\frac{1}{s}\right)\right]^\alpha\left(\int_0^s f_1^*(\tau)^{q_1}\,d\tau\right)^{\frac{1}{q_1}}\left(\int_0^s f_2^*(\tau)^{q_2}\,d\tau\right)^{\frac{1}{q_2}}\right\}\\
 & \le\sup_{s>0}\,\left\{\left[\log\left(e+\frac{1}{s}\right)\right]^{\alpha_1}\int_0^s f_1^*(\tau)^{q_1}\,d\tau\right\}^{\frac{1}{q_1}}\,\,
 \sup_{s>0}\,\left\{\left[\log\left(e+\frac{1}{s}\right)\right]^{\alpha_2}\int_0^s f_2^*(\tau)^{q_2}\,d\tau\right\}^{\frac{1}{q_2}}\\
 & =\|f_1\|_{{\mathfrak L}^{q_1,\infty}(\log{\mathfrak L})^{\alpha_1}}\|f_2\|_{{\mathfrak L}^{q_2,\infty}(\log{\mathfrak L})^{\alpha_2}}.
\end{align*}
Thus \eqref{eq:2.14} holds. 
Furthermore, for any $\tilde{f}_1\in {\mathfrak L}_{{\rm ul}}^{q_1,\infty}(\log{\mathfrak L})^{\alpha_1}$, $\tilde{f}_2\in {\mathfrak L}_{{\rm ul}}^{q_2,\infty}(\log{\mathfrak L})^{\alpha_2}$, and $\rho>0$, 
by \eqref{eq:2.14} we have 
\begin{align*}
|||\tilde{f}_1\tilde{f}_2|||_{1,\alpha;\rho}
 & =\sup_{x\in{\mathbb R}^n}\|\tilde{f}_1\tilde{f}_2\chi_{B(x,\rho)}\|_{{\mathfrak L}^{1,\infty}(\log{\mathfrak L})^\alpha}\\
 & \le \sup_{x\in{\mathbb R}^n}\left\{\|\tilde{f}_1\chi_{B(x,\rho)}\|_{{\mathfrak L}^{q_1,\infty}(\log{\mathfrak L})^{\alpha_1}}\|\tilde{f}_2\chi_{B(x,\rho)}\|_{{\mathfrak L}^{q_2,\infty}(\log{\mathfrak L})^{\alpha_2}}\right\}\\
 & \le \sup_{x\in{\mathbb R}^n}\|\tilde{f}_1\chi_{B(x,\rho)}\|_{{\mathfrak L}^{q_1,\infty}(\log{\mathfrak L})^{\alpha_1}}\cdot 
  \sup_{x\in{\mathbb R}^n}\|\tilde{f}_2\chi_{B(x,\rho)}\|_{{\mathfrak L}^{q_2,\infty}(\log{\mathfrak L})^{\alpha_2}}\\
 & =|||\tilde{f}_1|||_{q_1,\alpha_1;\rho}|||\tilde{f}_2|||_{q_2,\alpha_2;\rho}.
\end{align*}
Thus \eqref{eq:2.15} holds, and Lemma~\ref{Lemma:2.2} follows for $q_1$, $q_2\in[1,\infty)$.
If $q_1=\infty$ or $q_2=\infty$, the conclusion follows from \eqref{eq:1.2}. 
$\Box$
\begin{lemma}
\label{Lemma:2.3}
Let $q\in[1,\infty)$ and $\alpha\ge 0$. 
Then, for any $r>0$ with $rq\ge 1$, 
\begin{equation*}
\begin{array}{ll}
\||f|^r\|_{{\mathfrak L}^{q,\infty}(\log{\mathfrak L})^\alpha}=\|f\|_{{\mathfrak L}^{rq,\infty}(\log{\mathfrak L})^{\alpha}}^r
 & \quad\mbox{for $f\in {\mathfrak L}^{{rq},\infty}(\log{\mathfrak L})^{\alpha}$},\vspace{5pt}\\
||||\tilde{f}|^r|||_{q,\alpha;\rho}=|||\tilde{f}|||_{rq,\alpha;\rho}^r 
 & \quad\mbox{for $\tilde{f}\in {\mathfrak L}_{{\rm ul}}^{q,\infty}(\log{\mathfrak L})^{\alpha}$ and $\rho>0$}.
\end{array}
\end{equation*}
\end{lemma}
{\bf Proof.}
It follows from \eqref{eq:2.9} that 
\begin{align*}
\||f|^r\|_{{\mathfrak L}^{q,\infty}(\log{\mathfrak L})^\alpha} & 
=\sup_{s>0}\,\left\{\left[\log\left(e+\frac{1}{s}\right)\right]^\alpha \int_0^s ((|f|^r)^q)^*(\tau)\,d\tau\right\}^{\frac{1}{q}}\\
 & =\sup_{s>0}\,\left\{\left[\log\left(e+\frac{1}{s}\right)\right]^\alpha \int_0^s (|f|^{rq})^*(\tau)\,d\tau\right\}^{\frac{r}{rq}}
 =\|f\|_{{\mathfrak L}^{rq,\infty}(\log{\mathfrak L})^{\alpha}}^r
\end{align*}
for $f\in {\mathfrak L}^{{rq},\infty}(\log{\mathfrak L})^{\alpha}$.
Then 
\begin{align*}
||||\tilde{f}|^r|||_{q,\alpha;\rho}=\sup_{x\in{\mathbb R}^n}\| |\tilde{f}|^r\chi_{B(x,\rho)}\|_{{\mathfrak L}^{q,\infty}(\log{\mathfrak L})^\alpha}
=\sup_{x\in{\mathbb R}^n}\| |\tilde{f}|\chi_{B(x,\rho)}\|_{{\mathfrak L}^{rq,\infty}(\log{\mathfrak L})^{\alpha}}^r
=|||\tilde{f}|||_{rq,\alpha;\rho}^r
\end{align*}
for $\tilde{f}\in {\mathfrak L}_{{\rm ul}}^{{rq},\infty}(\log{\mathfrak L})^{\alpha}$ and $\rho>0$.
Thus Lemma~\ref{Lemma:2.3} follows. 
$\Box$
\subsection{Hardy's inequalities}
We recall the following two lemmas on Hardy's inequality.
(See \cite{Muckenhoup}*{Theorems~1 and 2}.)  
Throughout this paper, for any $q\in[1,\infty]$, we denote by $q'$ the H\"older conjugate of $q$, that is, 
$q'=q/(q-1)$ if $q\in(1,\infty)$, $q'=\infty$ if $q=1$, and $q'=1$ if $q=\infty$.
\begin{lemma}
\label{Lemma:2.4}
Let $q\in[1,\infty]$. 
Let $U$ and $V$ be locally integrable functions in $[0,\infty)$. 
Then there exists $C>0$ such that 
$$
\|U\tilde{F}\|_{L^q((0,\infty))}
\le C\|Vf\|_{L^q((0,\infty))}
\quad\mbox{with}\quad
\tilde{F}(s):=\int_0^s f(\tau)\,d\tau
$$
holds for all locally integrable function $f$ in $[0,\infty)$ 
if and only if 
$$
\sup_{s>0}\,
\left\{\|U\|_{L^q((s,\infty))}\|V^{-1}\|_{L^{q'}((0,s))}\right\}<\infty.
$$
\end{lemma}
\begin{lemma}
\label{Lemma:2.5}
Let $q\in[1,\infty]$. 
Let $U$ and $V$ be locally integrable functions in $[0,\infty)$. 
Then there exists $C>0$ such that 
$$
\|U\hat{F}\|_{L^q((0,\infty))}
\le C\|Vf\|_{L^q((0,\infty))}
\quad\mbox{with}\quad
\hat{F}(s):=\int_s^\infty f(\tau)\,d\tau
$$
holds for all locally integrable function $f$ in $(0,\infty)$ with $f\in L^1((1,\infty))$ if and only if 
$$
\sup_{s>0}\left\{\|U\|_{L^q((0,s))}\|V^{-1}\|_{L^{q'}((s,\infty))}\right\}<\infty.
$$
\end{lemma}
\subsection{Fundamental solutions}
Let $\theta\in(0,2]$. 
Let $G_\theta$ be the fundamental solution to the fractional heat equation
$$
\partial_t v+(-\Delta)^{\frac{\theta}{2}}v=0\quad\mbox{in}\quad {\mathbb R}^n\times(0,\infty).
$$
The function $G_\theta$ is positive and smooth in ${\mathbb R}^n\times(0,\infty)$ 
and it satisfies
\begin{equation}
\label{eq:2.16}
\begin{array}{ll}
\displaystyle{G_\theta(x,t)=(4\pi t)^{-\frac{n}{2}}\exp\left(-\frac{|x|^2}{4t}\right)\le Ch_{\theta,t}(x)} & \quad\mbox{if}\quad \theta=2,\vspace{5pt}\\
G_\theta(x,t)
\asymp h_{\theta,t}(x) & \quad\mbox{if}\quad 0<\theta<2,
\end{array}
\end{equation}
for $(x,t)\in{\mathbb R}^n\times(0,\infty)$, where
\begin{equation}
\label{eq:2.17}
h_{\theta,t}(x):=t^{-\frac{n}{\theta}}\left(1+t^{-\frac{1}{\theta}}|x|\right)^{-n-\theta}.
\end{equation}
Furthermore, 
\begin{eqnarray*}
 & \bullet & G_\theta(x,t)=t^{-\frac{n}{\theta}}G_\theta\left(t^{-\frac{1}{\theta}}x,1\right),
 \quad \int_{{\mathbb R}^n}G_\theta(x,t)\,dx=1,\\
 & \bullet &  \mbox{$G_\theta(\cdot,1)$ is radially symmetric and $G_\theta(x,1)\le G_\theta(y,1)$ if $|x|\ge |y|$},\\
 & \bullet & G_\theta(x,t)=\int_{{\mathbb R}^n}G_\theta(x-y,t-s)G_\theta(y,s)\,dy,
\end{eqnarray*}
for $x$, $y\in{\mathbb R}^n$ and $0<s<t$ 
(see e.g., \cites{BJ, BraK, S}) and 
\begin{equation}
\label{eq:2.18}
\lim_{t\to +0}\|S_\theta(t)\eta-\eta\|_{L^\infty}=0\quad\mbox{for}\quad \eta\in C_0({\mathbb R}^n).
\end{equation}
In addition, 
it follows from Young's inequality that 
\begin{equation}
\label{eq:2.19}
\|S_\theta(t)\eta\|_{L^{{q}}}\le Ct^{-\frac{n}{\theta}\left(\frac{1}{r}-\frac{1}{q}\right)}\|\eta\|_{L^{{r}}}
\end{equation}
for $\eta\in L^{{r}}$, $1\le r\le q\le\infty$, and $t>0$.
\section{Decay estimates of $S_\theta(t)\varphi$}
In this section we obtain decay estimates of $S_\theta(t)\varphi$ in our weak Zygmund type spaces 
and uniformly local weak Zygmund type spaces. 
For simplicity we write $g_t:=G_\theta(\cdot,t)$ and $h_t:=h_{\theta,t}$. 
\begin{proposition}
\label{Proposition:3.1}
Let $\theta\in(0,2]$, $1\le r\le q\le \infty$, and $\alpha$, $\beta\ge 0$. 
Assume that $\alpha\le\beta$ if $r=q$. 
Then there exists $C>0$ such that
$$
\|S_\theta(t)\varphi\|_{{\mathfrak L}^{q,\infty}(\log{\mathfrak L})^\beta}
\le Ct^{-\frac{n}{\theta}\left(\frac{1}{r}-\frac{1}{q}\right)}
\left[\log\left(e+\frac{1}{t}\right)\right]^{-\frac{\alpha}{r}+\frac{\beta}{q}}
\|\varphi\|_{{\mathfrak L}^{r,\infty}(\log{\mathfrak L})^\alpha}
$$
for $\varphi\in {\mathfrak L}^{r,\infty}(\log{\mathfrak L})^\alpha$ and $t>0$.
\end{proposition}
Before starting the proof, 
we recall the following relations on logarithmic functions: 
for any fixed $L>1$ and $k>0$, 
\begin{equation}
\label{eq:3.1}
\log\left(e+\frac{1}{s}\right)\asymp 
\log\left(L+\frac{1}{s}\right)\asymp
\log\left(e+\frac{k}{s}\right)\asymp
\log\left(e+\frac{1}{s^k}\right)
\quad\mbox{for}\quad s>0. 
\end{equation}
Furthermore, 
\begin{lemma}
\label{Lemma:3.1}
{\rm (1)} 
  Let $q>-1$ and $\alpha\in{\mathbb R}$. Then there exists $C_1>0$ such that 
  $$
  \int_0^s \tau^q\left[\log\left(e+\frac{1}{\tau}\right)\right]^\alpha\,d\tau\le C_1s^{q+1}\left[\log\left(e+\frac{1}{s}\right)\right]^\alpha
  \quad\mbox{for}\quad s>0.
  $$
{\rm (2)}
  Let $S>0$ and $\alpha<-1$. 
  Then there exists $C_2>0$ such that 
  $$
  \int_0^s \tau^{-1}\left[\log\left(e+\frac{1}{\tau}\right)\right]^\alpha\,d\tau\le C_2\left[\log\left(e+\frac{1}{s}\right)\right]^{\alpha+1}
  \quad\mbox{for}\quad s\in(0,S).
  $$
{\rm (3)}
  Let $q<-1$ and $\alpha\in{\mathbb R}$. Then there exists $C_3>0$ such that 
  $$
  \int_s^\infty \tau^q\left[\log\left(e+\frac{1}{\tau}\right)\right]^\alpha\,d\tau\le C_3s^{q+1}\left[\log\left(e+\frac{1}{s}\right)\right]^\alpha
  \quad\mbox{for}\quad s>0.
  $$
\end{lemma}
{\bf Proof.} 
We prove assertion~(1). 
Let $\delta>0$ be such that $q-\delta>-1$. 
Then there exists $L\in[e,\infty)$ such that 
\begin{equation}
\label{eq:3.2}
\mbox{the function $\displaystyle{(0,\infty)\ni \tau\mapsto \tau^\delta\left[\log\left(L+\frac{1}{\tau}\right)\right]^\alpha}$ is non-decreasing}.
\end{equation}
This together with \eqref{eq:3.1} implies that
\begin{align*}
\int_0^s \tau^q\left[\log\left(e+\frac{1}{\tau}\right)\right]^\alpha\,d\tau
 & \le C\int_0^s \tau^{q-\delta}\cdot\tau^\delta\left[\log\left(L+\frac{1}{\tau}\right)\right]^\alpha\,d\tau\\
 & \le Cs^\delta\left[\log\left(L+\frac{1}{s}\right)\right]^\alpha\int_0^s \tau^{q-\delta}\,d\tau
 \le Cs^{q+1}\left[\log\left(e+\frac{1}{s}\right)\right]^\alpha
\end{align*}
for $s>0$. Thus assertion~(1) follows. 

We prove assertion~(2). Let $S>0$. 
It follows that  
$$
\int_0^s \tau^{-1}\left[\log\left(e+\frac{1}{\tau}\right)\right]^\alpha\,d\tau
\le C\int_0^s \tau^{-1}|\log \tau|^\alpha \,d\tau\le C|\log s|^{\alpha+1}
\le C\left[\log\left(e+\frac{1}{s}\right)\right]^{\alpha+1}
$$
for $s\in(0,1/2)$. 
If $S\ge 1/2$, then
\begin{align*}
\int_0^s \tau^{-1}\left[\log\left(e+\frac{1}{\tau}\right)\right]^\alpha\,d\tau
&
\le \int_{\frac{1}{4}}^S \tau^{-1}\left[\log\left(e+\frac{1}{\tau}\right)\right]^\alpha\,d\tau+C
\\
&
\le C\le C\left[\log\left(e+\frac{1}{s}\right)\right]^{\alpha+1}
\end{align*}
for $s\in[1/2,S)$. Thus assertion~(2) follows. 

It remains to prove assertion~(3). 
Let $\delta'>0$ be such that $q+\delta'<-1$. 
Then there exists $L'\in[e,\infty)$ such that 
$$
\mbox{the function $\displaystyle{(0,\infty)\ni \tau\mapsto \tau^{-\delta'}\left[\log\left(L'+\frac{1}{\tau}\right)\right]^\alpha}$ is non-increasing}.
$$
This together with \eqref{eq:3.1} implies that 
\begin{align*}
\int_s^\infty \tau^q\left[\log\left(e+\frac{1}{\tau}\right)\right]^\alpha\,d\tau
 & \le C\int_s^\infty \tau^{q+\delta'}\cdot\tau^{-\delta'}\left[\log\left(L'+\frac{1}{\tau}\right)\right]^\alpha\,d\tau\\
 & \le Cs^{-\delta'}\left[\log\left(L'+\frac{1}{s}\right)\right]^\alpha\int_s^\infty \tau^{q+\delta'}\,d\tau
 \le Cs^{q+1}\left[\log\left(e+\frac{1}{s}\right)\right]^\alpha
\end{align*}
for $s>0$. Thus assertion~(3) follows. The proof is complete.
$\Box$\vspace{5pt}

Next, we prepare the following lemma on $h_t^{*}$, 
where $h_t=h_{\theta,t}$ is as in \eqref{eq:2.17}.

\begin{lemma}
\label{Lemma:3.2}
Let $1\le r\le q<\infty$ and $\gamma\in{\mathbb R}$. 
Assume that $\gamma\ge 0$ if $r=q$. 
Then there exists $C>0$ such that
\begin{equation}
\label{eq:3.3}
\int_0^\infty\tau^{q\left(1-\frac{1}{r}\right)}\left[\log\left(e+\frac{1}{\tau}\right)\right]^\gamma (h_t^*(\tau))^q\,d\tau
\le Ct^{-\frac{nq}{\theta}\left(\frac{1}{r}-\frac{1}{q}\right)}\left[\log\left(e+\frac{1}{t}\right)\right]^\gamma
\end{equation}
for $t>0$.
\end{lemma}
{\bf Proof.}
It follows from \eqref{eq:2.17} that
$$
(h_t)^*(s)=h_t\left(\left(\omega_n^{-1}s\right)^{\frac{1}{n}}e_1\right)
\le Ct^{-\frac{n}{\theta}}\left(1+t^{-\frac{1}{\theta}}s^{\frac{1}{n}}\right)^{-n-\theta}
$$
for $s\in[0,\infty)$ and $t\in(0,\infty)$, where $\omega_n$ is the volume of the $n$-dimensional unit ball $B(0,1)$ 
and $e_1:=(1,0,\dots,0)\in{\mathbb R}^n$.  
Then 
\begin{equation}
\label{eq:3.4}
\begin{split}
I & :=\int_0^\infty \tau^{q\left(1-\frac{1}{r}\right)}\left[\log\left(e+\frac{1}{\tau}\right)\right]^\gamma (h_t^*(\tau))^q\,d\tau\\
 & \,\le Ct^{-\frac{nq}{\theta}}\int_0^\infty \tau^{q\left(1-\frac{1}{r}\right)}\left[\log\left(e+\frac{1}{\tau}\right)\right]^\gamma
 \left(1+t^{-\frac{1}{\theta}}\tau^{\frac{1}{n}}\right)^{-q(n+\theta)}\,d\tau\\
 & \,\le Ct^{-\frac{nq}{\theta}\left(\frac{1}{r}-\frac{1}{q}\right)}
 \int_0^\infty \xi^{nq\left(1-\frac{1}{r}\right)+n-1}(1+\xi)^{-q(n+\theta)}
\left[\log\left(e+\frac{1}{(t^{1/\theta}\xi)^n}\right)\right]^\gamma\,d\xi
\end{split}
\end{equation}
for $t>0$. 

We first consider the case of $\gamma\ge 0$. 
It follows from \eqref{eq:3.1} that
\begin{equation}
\label{eq:3.5}
\begin{split}
\left[\log\left(e+\frac{1}{(t^{1/\theta}\xi)^n}\right)\right]^\gamma
 & \le C\left[\log\left(e+\frac{1}{t^{1/\theta}\xi}\right)\right]^\gamma\\
 & \le C\left[\log\left(e+\frac{1}{t^{1/\theta}}\right)+\log\left(e+\frac{1}{\xi}\right)\right]^\gamma\\
 & \le C\left[\log\left(e+\frac{1}{t}\right)\right]^\gamma
 +C\left[\log\left(e+\frac{1}{\xi}\right)\right]^\gamma
\end{split}
\end{equation}
for $t>0$ and $\xi\in(0,1/2)$. 
Similarly, by \eqref{eq:3.1} we have
\begin{equation}
\label{eq:3.6}
\left[\log\left(e+\frac{1}{(t^{1/\theta}\xi)^n}\right)\right]^\gamma
\le \left[\log\left(e+\frac{2^n}{(t^{1/\theta})^n}\right)\right]^\gamma
 \le C\left[\log\left(e+\frac{1}{t}\right)\right]^\gamma
\end{equation}
for $t>0$ and $\xi\in[1/2,\infty)$. 
Since 
\begin{equation}
\label{eq:3.7}
nq\left(1-\frac{1}{r}\right)+n-1-q(n+\theta)=-\frac{nq}{r}+n-1-q\theta
=-nq\left(\frac{1}{r}-\frac{1}{q}\right)-1-q\theta<-1,
\end{equation}
by Lemma~\ref{Lemma:3.1}, \eqref{eq:3.4}, \eqref{eq:3.5}, and \eqref{eq:3.6} we obtain 
\begin{align*}
I & \le Ct^{-\frac{nq}{\theta}\left(\frac{1}{r}-\frac{1}{q}\right)}
\int_0^{1/2} \xi^{nq\left(1-\frac{1}{r}\right)+n-1}
\left(\left[\log\left(e+\frac{1}{t}\right)\right]^\gamma
+\left[\log\left(e+\frac{1}{\xi}\right)\right]^\gamma\right)\,d\xi\\
& \qquad\quad
+Ct^{-\frac{nq}{\theta}\left(\frac{1}{r}-\frac{1}{q}\right)}\int_{1/2}^\infty \xi^{nq\left(1-\frac{1}{r}\right)+n-1}(1+\xi)^{-q(n+\theta)}\left[\log\left(e+\frac{1}{t}\right)\right]^\gamma\,d\xi\\
& \le Ct^{-\frac{nq}{\theta}\left(\frac{1}{r}-\frac{1}{q}\right)}\left(1+\left[\log\left(e+\frac{1}{t}\right)\right]^\gamma\right)\le Ct^{-\frac{nq}{\theta}\left(\frac{1}{r}-\frac{1}{q}\right)}\left[\log\left(e+\frac{1}{t}\right)\right]^\gamma
\end{align*}
for $t>0$. This implies \eqref{eq:3.3} in the case of $\gamma\ge 0$. 

Consider the case of $\gamma<0$. 
Then, by \eqref{eq:3.1} we have
\begin{equation}
\label{eq:3.8}
\left[\log\left(e+\frac{1}{(t^{1/\theta}\xi)^n}\right)\right]^\gamma
\le \left[\log\left(e+\frac{2^n}{(t^{1/\theta})^n}\right)\right]^\gamma
 \le C \left[\log\left(e+\frac{1}{t}\right)\right]^\gamma
\end{equation}
for $t>0$ and $\xi\in(0,1/2)$. 
Let $0<\delta<\theta q/|\gamma|$.
We find $L\in[e,\infty)$ such that 
the function $f$ in $(0,\infty)$ defined by 
$$
f(z):=z^\delta\log\left(L+\frac{1}{z^n}\right)
$$ 
is non-decreasing in $(0,\infty)$.
Since $\gamma<0$, by \eqref{eq:3.1} we obtain
\begin{align*}
\left[\log\left(e+\frac{1}{(t^{1/\theta}\xi)^n}\right)\right]^\gamma & 
\le C\left[\log\left(L+\frac{1}{(t^{1/\theta}\xi)^n}\right)\right]^\gamma
=C\left[z^{-\delta\gamma}f(z)^\gamma\right]\biggr|_{z=t^{\frac{1}{\theta}}\xi}\\
 & \le C(t^{\frac{1}{\theta}}\xi)^{-\delta \gamma}f(z)^\gamma\biggr|_{z=t^{\frac{1}{\theta}}/2}
 \le C\xi^{-\delta\gamma}\left[\log\left(e+\frac{2^n}{t^{n/\theta}}\right)\right]^\gamma\\
 & \le C\xi^{-\delta\gamma}\left[\log\left(e+\frac{1}{t}\right)\right]^\gamma
\end{align*}
for $t>0$ and $\xi\in[1/2,\infty)$. 
This together with \eqref{eq:3.7} and \eqref{eq:3.8} implies that 
\begin{align*}
I & \le Ct^{-\frac{nq}{\theta}\left(\frac{1}{r}-\frac{1}{q}\right)}
\int_0^{1/2} \xi^{nq\left(1-\frac{1}{r}\right)+n-1}\left[\log\left(e+\frac{1}{t}\right)\right]^\gamma\,d\xi\\
 & \qquad\quad
 +Ct^{-\frac{nq}{\theta}\left(\frac{1}{r}-\frac{1}{q}\right)}
\int_{1/2}^\infty \xi^{nq\left(1-\frac{1}{r}\right)+n-1-\delta\gamma}(1+\xi)^{-q(n+\theta)}\left[\log\left(e+\frac{1}{t}\right)\right]^\gamma\,d\xi\\
 & \le Ct^{-\frac{nq}{2}\left(\frac{1}{r}-\frac{1}{q}\right)}\left[\log\left(e+\frac{1}{t}\right)\right]^\gamma
\end{align*}
for $t>0$. This implies \eqref{eq:3.3} in the case of $\gamma<0$.
Thus Lemma~\ref{Lemma:3.2} follows.
$\Box$
\vspace{5pt}
\newline
{\bf Proof of Proposition~\ref{Proposition:3.1}.}
The proof is divided into the following three cases: 
$$
1\le r<q<\infty,\qquad 1\le r=q<\infty,\qquad 1\le r\le q=\infty.
$$ 
\underline{Step 1.}
Consider the case of $1\le r<q<\infty$. 
It follows from \eqref{eq:2.4}, \eqref{eq:2.7}, \eqref{eq:2.9}, and \eqref{eq:2.16} that 
\begin{align*}
\left\|S_\theta(t)\varphi\right\|_{{\mathfrak L}^{q,\infty}(\log{\mathfrak L})^\beta}^q
 & =\sup_{s>0}\,\left\{\left[\log\left(e+\frac{1}{s}\right)\right]^{\beta}\int_0^s \left(\left(S_\theta(t)\varphi\right)^*(\tau)\right)^q\,d\tau\right\}\\
 & \le \sup_{s>0}\,\left\{\left[\log\left(e+\frac{1}{s}\right)\right]^{\beta}\int_0^s \left(\left(S_\theta(t)\varphi\right)^{**}(\tau)\right)^q\,d\tau\right\}\\
 & \le \sup_{s>0}\,\left\{\left[\log\left(e+\frac{1}{s}\right)\right]^{\beta}\int_0^s \left(\int_\tau^\infty g_t^{**}(\eta)\varphi^{**}(\eta)\,d\eta\right)^q\,d\tau\right\}\\
 & \le C\sup_{s>0}\,\left\{\left[\log\left(e+\frac{1}{s}\right)\right]^{\beta}\int_0^s \left(\int_\tau^\infty h_t^{**}(\eta)\varphi^{**}(\eta)\,d\eta\right)^q\,d\tau\right\}
\end{align*}
for $t>0$. Furthermore, thanks to \eqref{eq:2.11}, we have
\begin{equation}
\label{eq:3.9}
\left\|S_\theta(t)\varphi\right\|_{{\mathfrak L}^{q,\infty}(\log{\mathfrak L})^\beta}^q
\le C\int_0^\infty \left(\left[\log\left(e+\frac{1}{\tau}\right)\right]^{\frac{\beta}{q}}\int_\tau^\infty h_t^{**}(\eta)\varphi^{**}(\eta)\,d\eta\right)^q\,d\tau
\end{equation}
for $t>0$. 
On the other hand, set
$$
U(\tau):=\left[\log\left(e+\frac{1}{\tau}\right)\right]^{\frac{\beta}{q}},
\quad
V(\tau):=\tau\left[\log\left(e+\frac{1}{\tau}\right)\right]^{\frac{\beta}{q}}, 
$$
for $\tau>0$. It follows from Lemma~\ref{Lemma:3.1}~(1), (3) that
\begin{align*}
 & \sup_{s>0}\,\left(\int_0^s |U(\tau)|^q\,d\tau\right)^{\frac{1}{q}}\left(\int_s^\infty |V(\tau)|^{-q'}\,d\tau\right)^{\frac{1}{q'}}\\
 & \le \sup_{s>0}\,\left\{Cs^{\frac{1}{q}}\left[\log\left(e+\frac{1}{s}\right)\right]^{\frac{\beta}{q}}\cdot 
 Cs^{-1+\frac{1}{q'}}\left[\log\left(e+\frac{1}{s}\right)\right]^{-\frac{\beta}{q}}\right\}<\infty.
\end{align*}
Then, by Lemma~\ref{Lemma:2.5}, \eqref{eq:2.3}, and \eqref{eq:3.9} we have 
$$
\begin{aligned}
\left\|S_\theta(t)\varphi\right\|_{{\mathfrak L}^{q,\infty}(\log{\mathfrak L})^\beta}^q
 & \le C\int_0^\infty \left(\tau\left[\log\left(e+\frac{1}{\tau}\right)\right]^{\frac{\beta}{q}} h_t^{**}(\tau)\varphi^{**}(\tau)\right)^q\,d\tau\\
 & \le C\sup_{s>0}\,\left\{\left[\log\left(e+\frac{1}{s}\right)\right]^{\alpha}s(\varphi^{**}(s))^r\right\}^{\frac{q}{r}}\\
 & \qquad\quad
 \times\int_0^\infty\bigg( \tau^{1-\frac{1}{r}}\left[\log\left(e+\frac{1}{\tau}\right)\right]^{-\frac{\alpha}{r}+\frac{\beta}{q}}h_t^{**}(\tau)\bigg)^q\,d\tau
\end{aligned}
$$
for $t>0$.
This together with \eqref{eq:2.5} and \eqref{eq:2.9} implies that
\begin{equation}
\label{eq:3.10}
\left\|S_\theta(t)\varphi\right\|_{{\mathfrak L}^{q,\infty}(\log{\mathfrak L})^\beta}^q
\le C\|\varphi\|_{{\mathfrak L}^{r,\infty}(\log{\mathfrak L})^\alpha}^q
  \int_0^\infty\left(\tau^{-\frac{1}{r}}\left[\log\left(e+\frac{1}{\tau}\right)\right]^{{\gamma}} \int_0^\tau h_t^*(\xi)\,d\xi\right)^q\,d\tau
\end{equation}
for $t>0$,
where 
$$
{\gamma}:=-\frac{\alpha}{r}+\frac{\beta}{q}.
$$
Set
$$
\tilde{U}(\tau)=\tau^{-\frac{1}{r}}\left[\log\left(e+\frac{1}{\tau}\right)\right]^{{\gamma}},
\quad
\tilde{V}(\tau)=\tau^{1-\frac{1}{r}}\left[\log\left(e+\frac{1}{\tau}\right)\right]^{{\gamma}}.
$$
Since $q>r$ and $q'<r'$, 
by Lemma~\ref{Lemma:3.1}~(1), (3) we have 
\begin{equation}
\label{eq:3.11}
\begin{split}
 & \sup_{s>0}\,\left(\int_s^\infty |\tilde{U}(\tau)|^q\,d\tau\right)^{\frac{1}{q}}
\left(\int_0^s |\tilde{V}(\tau)|^{-q'}\,{d\tau}\right)^{\frac{1}{q'}}\\
 & =\sup_{s>0}\,\left(\int_s^\infty \tau^{-\frac{q}{r}}\left[\log\left(e+\frac{1}{\tau}\right)\right]^{q{\gamma}}\,d\tau\right)^{\frac{1}{q}}
\left(\int_0^s \tau^{-\frac{q'}{r'}}\left[\log\left(e+\frac{1}{\tau}\right)\right]^{-q'{\gamma}}\,{d\tau}\right)^{\frac{1}{q'}}\\
& \le \sup_{s>0}\,\left\{C s^{\frac{1}{q}-\frac{1}{r}}\left[\log\left(e+\frac{1}{s}\right)\right]^{{\gamma}}\cdot Cs^{\frac{1}{q'}-\frac{1}{r'}}\left[\log\left(e+\frac{1}{s}\right)\right]^{-{\gamma}}\right\}<\infty.
\end{split}
\end{equation}
Applying Lemma~\ref{Lemma:2.4} to \eqref{eq:3.10}, by \eqref{eq:3.11}
we obtain 
\begin{equation*}
\begin{split}
 & \left\|S_\theta(t)\varphi\right\|_{{\mathfrak L}^{q,\infty}(\log{\mathfrak L})^\beta}^q\\
 & \le C\|\varphi\|_{{\mathfrak L}^{r,\infty}(\log{\mathfrak L})^\alpha}^q
  \int_0^\infty\left(\tau^{1-\frac{1}{r}}\left[\log\left(e+\frac{1}{\tau}\right)\right]^{\gamma}h_t^*(\tau)\right)^q\,d\tau\\
\end{split}
\end{equation*}
for $t>0$. 
This together with Lemma~\ref{Lemma:3.2} implies that
$$
\left\|S_\theta(t)\varphi\right\|_{{\mathfrak L}^{q,\infty}(\log{\mathfrak L})^\beta}^q
\le Ct^{-\frac{nq}{\theta}\left(\frac{1}{r}-\frac{1}{q}\right)}\left[\log\left(e+\frac{1}{t}\right)\right]^{q\gamma}\|\varphi\|_{{\mathfrak L}^{r,\infty}(\log{\mathfrak L})^\alpha}^q
$$
for $t>0$. Thus Proposition~\ref{Proposition:3.1} follows in the case of $1\le r<q<\infty$. 
\newline
\underline{Step 2.}
Consider the case of $1\le r=q<\infty$.  
It follows from H\"older's inequality and \eqref{eq:2.16} that 
\begin{align*}
\left|\left[S_\theta(t)\varphi\right](x)\right|^r
 & \le C\left(\int_{{\mathbb R}^n} |h_t(x-y)||\varphi(y)|\,dy\right)^r\\
 & \le C\left(\int_{{\mathbb R}^n}|h_t(x-y)|\,dy\right)^{r-1}\int_{{\mathbb R}^n} |h_t(x-y)||\varphi(y)|^r\,dy\\
 & \le C\int_{{\mathbb R}^n} |h_t(x-y)||\varphi(y)|^r\,dy.
\end{align*}
Then it follows from \eqref{eq:2.7} and \eqref{eq:2.9} that 
\begin{equation}
\label{eq:3.12}
\begin{split}
\|S_\theta(t)\varphi\|_{{\mathfrak L}^{r,\infty}(\log{\mathfrak L})^\beta}^r
 & =\sup_{s>0}\, \left\{\left[\log\left(e+\frac{1}{s}\right)\right]^{\beta} s (|S_\theta(t)\varphi|^r)^{**}(s)\right\}\\
 & \le \sup_{s>0}\, \left\{\left[\log\left(e+\frac{1}{s}\right)\right]^{\beta}s\int_s^\infty (h_t)^{**}(\tau)(|\varphi|^r)^{**}(s)\,d\tau\right\}
\end{split}
\end{equation}
for $t>0$. Set 
$$
\hat{U}(r)=r\left[\log\left(e+\frac{1}{r}\right)\right]^{\beta},\quad \hat{V}(r)=r^2\left[\log\left(e+\frac{1}{r}\right)\right]^{\beta}. 
$$
Similarly to \eqref{eq:3.2}, we find $L\in[e,\infty)$ such that
$$
\mbox{the function $\displaystyle{(0,\infty)\ni r\mapsto r\left[\log\left(L+\frac{1}{r}\right)\right]^{\beta}}$ is non-decreasing}.
$$
Then, by \eqref{eq:3.1} we have 
$$
\|\hat{U}\|_{L^\infty(0,s)} 
\le C\sup_{r\in(0,s)}\left\{r\left[\log\left(L+\frac{1}{r}\right)\right]^{\beta}\right\} \le Cs\left[\log\left(L+\frac{1}{s}\right)\right]^{\beta}
\le Cs\left[\log\left(e+\frac{1}{s}\right)\right]^{\beta}
$$
for $s>0$. This together with Lemma~\ref{Lemma:3.1}~(3) implies that 
\begin{equation}
\label{eq:3.13}
\begin{split}
 & \sup_{s>0}\,\left\{\|\hat{U}\|_{L^\infty((0,s))}\int_s^\infty |\hat{V}(\tau)|^{-1}\,d\tau\right\}\\
 & \le \sup_{s>0}\left\{Cs\left[\log\left(e+\frac{1}{s}\right)\right]^{\beta}\cdot Cs^{-1}\left[\log\left(e+\frac{1}{s}\right)\right]^{-\beta}\right\}<\infty.
\end{split}
\end{equation}
Applying Lemma~\ref{Lemma:2.5} with $q=\infty$, by \eqref{eq:2.9}, \eqref{eq:3.12} and \eqref{eq:3.13} we obtain 
\begin{equation}
\label{eq:3.14}
\begin{split}
 & \|S_\theta(t)\varphi\|_{{\mathfrak L}^{r,\infty}(\log{\mathfrak L})^\beta}^r\\
 & \le C\sup_{s>0}\, \left\{s^2\left[\log\left(e+\frac{1}{s}\right)\right]^{\beta} (h_t)^{**}(s)(|\varphi|^r)^{**}(s)\right\}\\
 & \le C\sup_{s>0}\, \left\{s\left[\log\left(e+\frac{1}{s}\right)\right]^{\beta-\alpha}(h_t)^{**}(s)\right\}
\cdot
\sup_{s>0}\left\{\,s\left[\log\left(e+\frac{1}{s}\right)\right]^{\alpha}(|\varphi|^r)^{**}(s)\right\}\\
& =C\|h_t\|_{{\mathfrak L}^{1,\infty}(\log{\mathfrak L})^{\beta-\alpha}}\|\varphi\|^{{r}}_{{{\mathfrak L}^{r,\infty}(\log{\mathfrak L})^\alpha}}
\end{split}
\end{equation}
for $t>0$. 
Furthermore, 
since $\alpha\le\beta$,
by Lemma~\ref{Lemma:3.2}, \eqref{eq:2.9}, and \eqref{eq:2.11}
we have
\begin{align*}
\|h_t\|_{{\mathfrak L}^{1,\infty}(\log{\mathfrak L})^{\beta-\alpha}}
 & =\sup_{s>0}\, \left\{\left[\log\left(e+\frac{1}{s}\right)\right]^{\beta-\alpha}\int_0^s (h_t)^*(\tau)\,d\tau\right\}\\
 & \le\int_0^\infty \left[\log\left(e+\frac{1}{\tau}\right)\right]^{\beta-\alpha}(h_t)^*(\tau)\,d\tau
 \le C\left[\log\left(e+\frac{1}{t}\right)\right]^{\beta-\alpha}
\end{align*}
for $t>0$. 
This together with \eqref{eq:3.14} implies that 
$$
\|S_\theta(t)\varphi\|_{{\mathfrak L}^{r,\infty}(\log{\mathfrak L})^\beta}^r
\le C\left[\log\left(e+\frac{1}{t}\right)\right]^{\beta-\alpha}\|\varphi\|^r_{{\mathfrak L}^{r,\infty}(\log{\mathfrak L})^\alpha}
\quad\mbox{for}\quad t>0. 
$$
Thus Proposition~\ref{Proposition:3.1} follows in the case of $1\le r=q<\infty$. 
\newline
\underline{Step 3.}
It remains to consider the case of $1\le r\le q=\infty$. 
If $r=q=\infty$, then it follows from \eqref{eq:2.16} that
$$
\|S_\theta(t)\varphi\|_{L^\infty}\le C\|\varphi\|_{L^\infty}\int_{{\mathbb R}^n}h_t(y)\,dy
\le C\|\varphi\|_{L^\infty}
$$
for $t>0$, and Proposition~\ref{Proposition:3.1} follows. 
On the other hand, in the case of $1\le r<q=\infty$, let $\tilde{q}\in(r,\infty)$. 
Then, by Proposition~\ref{Proposition:3.1} with $q=\tilde{q}>r$, \eqref{eq:1.3}, \eqref{eq:2.19}, and \eqref{eq:3.1} we have 
\begin{align*}
\|S_\theta(t)\varphi\|_{L^\infty} & =\left\|S_\theta\left(\frac{t}{2}\right)S_\theta\left(\frac{t}{2}\right)\varphi\right\|_{L^\infty}\\
 & \le Ct^{-\frac{n}{\theta \tilde{q}}}\left\|S_\theta\left(\frac{t}{2}\right)\varphi\right\|_{L^{\tilde{q}}}
 =C^{-\frac{n}{\theta \tilde{q}}}\left\|S_\theta\left(\frac{t}{2}\right)\varphi\right\|_{{\mathfrak L}^{\tilde{q},\infty}(\log{\mathfrak L})^0}\\
 & \le Ct^{-\frac{n}{\theta \tilde{q}}}\cdot Ct^{-\frac{n}{\theta}\left(\frac{1}{r}-\frac{1}{\tilde{q}}\right)}
 \left[\log\left(e+\frac{2}{t}\right)\right]^{-\frac{\alpha}{r}}\|\varphi\|_{{\mathfrak L}^{r,\infty}(\log{\mathfrak L})^\alpha}\\
 & =Ct^{-\frac{n}{\theta r}}\left[\log\left(e+\frac{1}{t}\right)\right]^{-\frac{\alpha}{r}}\|\varphi\|_{{\mathfrak L}^{r,\infty}(\log{\mathfrak L})^\alpha}
\end{align*}
for $t>0$. Thus Proposition~\ref{Proposition:3.1} follows in the case of $1\le r<q=\infty$. 
The proof of Proposition~\ref{Proposition:3.1} is complete.
$\Box$\vspace{5pt}

Furthermore, by Proposition~\ref{Proposition:3.1} we employ the arguments in the proof of \cite{HI01}*{Lemma~2.1} 
to obtain decay estimates of $S_\theta(t)\varphi$ in uniformly local weak Zygmund type spaces.  
\begin{proposition}
\label{Proposition:3.2}
Let $\theta\in(0,2]$, $1\le r\le q\le \infty$, and $\alpha$, $\beta\ge 0$. 
Assume that $\alpha\le\beta$ if $r=q$. 
There exists $C>0$ such that, for any $T>0$, 
\begin{equation}
\label{eq:3.15}
|||S_\theta(t)\varphi|||_{q,\beta;T^{1/\theta}}
\le Ct^{-\frac{n}{\theta}\left(\frac{1}{r}-\frac{1}{q}\right)}\left[\log\left(e+\frac{1}{t}\right)\right]^{-\frac{\alpha}{r}+\frac{\beta}{q}}|||\varphi|||_{r,\alpha;T^{1/\theta}}
\end{equation}
for $\varphi\in {\mathfrak L}^{r,\infty}_{{\rm ul}}(\log{\mathfrak L})^\alpha$ and $t\in(0,T]$.
\end{proposition}
{\bf Proof.}
We first consider the case of $\theta\in(0,2)$. 
It suffices to prove 
\begin{equation}
\label{eq:3.16}
t^{\frac{n}{\theta}\left(\frac{1}{r}-\frac{1}{q}\right)}\left[\log\left(e+\frac{1}{t}\right)\right]^{\frac{\alpha}{r}-\frac{\beta}{q}}
\left\|\chi_{B(z,T^{1/\theta})}S_\theta(t)\varphi\right\|_{{\mathfrak L}^{q,\infty}(\log{\mathfrak L})^\beta}\le C|||\varphi|||_{r,\alpha;T^{1/\theta}}
\end{equation}
for $z\in{\mathbb R}^n$ and $0<t\le T$. 
For the proof, by translating if necessary, we have only to consider the case of $z=0$. 

By Besicovitch's covering lemma 
we can find an integer $m$ depending only on $n$ and 
a set $\{x_{k,i}\}_{k=1,\dots,m,\,i\in{\mathbb N}}\subset{\mathbb R}^n\setminus B(0,10T^{1/\theta})$ such that 
\begin{equation}
\label{eq:3.17}
B_{k,i}\cap B_{k,j}=\emptyset\quad\mbox{if $i\not=j$}
\qquad\mbox{and}\qquad
{\mathbb R}^n\setminus B(0,10T^{1/\theta})\subset\bigcup_{k=1}^m\bigcup_{i=1}^\infty B_{k,i},
\end{equation}
where $B_{k,i}:=\overline{B(x_{k,i},T^{1/\theta})}$. 
Then 
\begin{equation}
\label{eq:3.18}
\left|\left[S_\theta(t)\varphi\right](x)\right|\le |u_0(x,t)|
+\sum_{k=1}^m\sum_{i=1}^\infty |u_{k,i}(x,t)|,
\quad(x,t)\in{\mathbb R}^n\times(0,T),
\end{equation}
where 
$$
u_0(x,t):=[S_\theta(t)(\varphi\chi_{B(0,10T^{1/\theta})})](x),
\quad
u_{k,i}(x,t):=[S_\theta(t)(\varphi\chi_{B_{k,i}})](x).
$$
By Proposition~\ref{Proposition:3.1} and \eqref{eq:1.4} we have 
\begin{equation}
\label{eq:3.19}
\begin{split}
 & \|u_0(t)\chi_{B(0,T^{1/\theta})}\|_{{\mathfrak L}^{q,\infty}(\log{\mathfrak L})^\beta}\le\|u_0(t)\|_{{\mathfrak L}^{q,\infty}(\log{\mathfrak L})^\beta}\\
 & \qquad\quad
 \le Ct^{-\frac{n}{\theta}\left(\frac{1}{r}-\frac{1}{q}\right)}
 \left[\log\left(e+\frac{1}{t}\right)\right]^{-\frac{\alpha}{r}+\frac{\beta}{q}}\left\|\varphi\chi_{B(0,10T^{1/\theta})}\right\|_{{\mathfrak L}^{r,\infty}(\log{\mathfrak L})^\alpha}\\
 & \qquad\quad
 \le Ct^{-\frac{n}{\theta}\left(\frac{1}{r}-\frac{1}{q}\right)}\left[\log\left(e+\frac{1}{t}\right)\right]^{-\frac{\alpha}{r}+\frac{\beta}{q}}|||\varphi|||_{r,\alpha;10T^{1/\theta}}\\
 & \qquad\quad
 \le Ct^{-\frac{n}{\theta}\left(\frac{1}{r}-\frac{1}{q}\right)}\left[\log\left(e+\frac{1}{t}\right)\right]^{-\frac{\alpha}{r}+\frac{\beta}{q}}|||\varphi|||_{r,\alpha;T^{1/\theta}}
\end{split}
\end{equation}
for $t\in(0,T]$. 

Let $k=1,\dots,m$ and $i{\in \mathbb{N}}$. 
By \eqref{eq:2.16} we have 
\begin{equation}
\label{eq:3.20}
\begin{split}
|u_{k,i}(x,t)| & \le C\int_{B(x_{k,i},T^{1/\theta})}h_t(x-y)|\varphi(y)|\,dy\\
 & =C\int_{\mathbb R^n}h_t(x-z-x_{k,i})\varphi_{k,i}(z)\,dz
\end{split}
\end{equation}
for $(x,t)\in{\mathbb R}^n\times(0,\infty)$, where 
$\varphi_{k,i}(x)=|\varphi(x+x_{k,i})|\chi_{B(0,T^{1/\theta})}$. 
Since $|x_{k,i}|\ge 10T^{1/\theta}$, it follows that 
\begin{align*}
 & (1+T^{-\frac{1}{\theta}}|x_{k,i}|)(1+t^{-\frac{1}{\theta}}|x-z|)\\
 & =1+T^{-\frac{1}{\theta}}|x_{k,i}|+t^{-\frac{1}{\theta}}|x-z|+t^{-\frac{1}{\theta}}T^{-\frac{1}{\theta}}|x_{k,i}||x-z|\\
 & \le 1+3t^{-\frac{1}{\theta}}|x_{k,i}|+t^{-\frac{1}{\theta}}|x-z|\\
 & = 1+4t^{-\frac{1}{\theta}}(|x_{k,i}|-|x-z|)-t^{-\frac{1}{\theta}}|x_{k,i}|+5t^{-\frac{1}{\theta}}|x-z|\\
 & \le  4(1+t^{-\frac{1}{\theta}}(|x_{k,i}|-|x-z|))\le 4(1+t^{-\frac{1}{\theta}}|x-z-x_{k,i}|)
\end{align*}
for $x$, $z\in B(0,T^{1/\theta})$ and $t\in(0,T)$. 
This together with \eqref{eq:2.16} implies that 
\begin{equation}
\label{eq:3.21}
\begin{split}
h_t(x-z-x_{k,i}) & \le Ct^{-\frac{n}{\theta}}(1+T^{-\frac{1}{\theta}}|x_{k,i}|)^{-n-\theta}(1+t^{-\frac{1}{\theta}}|x-z|)^{-n-\theta}\\
 & \le C(1+T^{-\frac{1}{\theta}}|x_{k,i}|)^{-n-\theta}g_t(x-z)
\end{split}
\end{equation}
for $x$, $z\in B(0,T^{1/\theta})$ and $t\in(0,T)$. 
We observe from \eqref{eq:3.20} and \eqref{eq:3.21} that 
$$
|u_{k,i}(x,t)|\le C(1+T^{-\frac{1}{\theta}}|x_{k,i}|)^{-n-\theta}
[S_\theta(t)\varphi_{k,i}](x)
$$
for $x\in B(0,T^{1/\theta})$ and $t\in(0,T)$. 
Then, by Proposition~\ref{Proposition:3.1} we obtain 
\begin{equation}
\label{eq:3.22}
\begin{split}
 & \|u_{k,i}(t)\chi_{B(0,T^{1/\theta})}\|_{{\mathfrak L}^{q,\infty}(\log{\mathfrak L})^\beta}\\
 & \le C(1+T^{-\frac{1}{\theta}}|x_{k,i}|)^{-n-\theta}\|S_\theta(t)\varphi_{k,i}\|_{{\mathfrak L}^{q,\infty}(\log{\mathfrak L})^\beta}\\
 & \le C(1+T^{-\frac{1}{\theta}}|x_{k,i}|)^{-n-\theta}t^{-\frac{n}{\theta}\left(\frac{1}{r}-\frac{1}{q}\right)}\left[\log\left(e+\frac{1}{t}\right)\right]^{-\frac{\alpha}{r}+\frac{\beta}{q}}
 \|\varphi_{k,i}\|_{{\mathfrak L}^{r,\infty}(\log{\mathfrak L})^\alpha}\\
 & =C(1+T^{-\frac{1}{\theta}}|x_{k,i}|)^{-n-\theta}t^{-\frac{n}{\theta}\left(\frac{1}{r}-\frac{1}{q}\right)}\left[\log\left(e+\frac{1}{t}\right)\right]^{-\frac{\alpha}{r}+\frac{\beta}{q}}
 \|\varphi\chi_{B(x_{k,i},T^{1/\theta})}\|_{{\mathfrak L}^{r,\infty}(\log{\mathfrak L})^\alpha}\\
 & \le C(1+T^{-\frac{1}{\theta}}|x_{k,i}|)^{-n-\theta}t^{-\frac{n}{\theta}\left(\frac{1}{r}-\frac{1}{q}\right)}\left[\log\left(e+\frac{1}{t}\right)\right]^{-\frac{\alpha}{r}+\frac{\beta}{q}}
 |||\varphi|||_{r,\alpha;T^{1/\theta}}
\end{split}
\end{equation}
for $t\in(0,T)$. 

On the other hand, since
$$
\frac{|y|}{2}\le\frac{1}{2}\left(|x_{k,i}|+T^{\frac{1}{\theta}}\right)\le|x_{k,i}|\quad\mbox{for}\quad y\in B_{k,i},
$$
we have 
$$
\frac{1}{|B_{k,i}|}\int_{B_{k,i}} \left(1+\frac{1}{2}T^{-\frac{1}{\theta}}|y|\right)^{-n-\theta}\,dy\ge (1+T^{-\frac{1}{\theta}}|x_{k,i}|)^{-n-\theta}.
$$
Then, by \eqref{eq:3.17} we see that 
\begin{equation}
\label{eq:3.23}
\begin{split}
\sum_{i=1}^\infty (1+T^{-\frac{1}{\theta}}|x_{k,i}|)^{-n-\theta}
 & \le CT^{-\frac{n}{\theta}}\sum_{i=1}^\infty\int_{B_{k,i}} \left(1+\frac{1}{2}T^{-\frac{1}{\theta}}|y|\right)^{-n-\theta}\,dy\\
 & \le CT^{-\frac{n}{\theta}}\int_{{\mathbb R}^n} \left(1+\frac{1}{2}T^{-\frac{1}{\theta}}|y|\right)^{-n-\theta}\,dy
\le C
\end{split}
\end{equation}
for $T>0$. 
Combining \eqref{eq:3.18}, \eqref{eq:3.19}, \eqref{eq:3.22}, and \eqref{eq:3.23}
we obtain 
\begin{align*}
 & t^{\frac{n}{\theta}\left(\frac{1}{r}-\frac{1}{q}\right)}
 \left[\log\left(e+\frac{1}{t}\right)\right]^{\frac{\alpha}{r}-\frac{\beta}{q}}\left\|\chi_{B(0,T^{1/\theta})}S_\theta(t)\varphi\right\|_{{\mathfrak L}^{q,\infty}(\log{\mathfrak L})^\beta}\\
 & \le C|||\varphi|||_{r,\alpha;T^{1/\theta}}+C|||\varphi|||_{r,\alpha;T^{1/\theta}}
 \sum_{k=1}^m\sum_{i=1}^\infty (1+T^{-\frac{1}{\theta}}|x_{k,i}|)^{-n-\theta}
 \le C|||\varphi|||_{r,\alpha;T^{1/\theta}}
\end{align*}
for $t\in(0,T)$. 
This implies \eqref{eq:3.16} with $z=0$, that is, \eqref{eq:3.15} holds. 
Thus Proposition~\ref{Proposition:3.2} follows in the case of $0<\theta<2$. 

Consider the case of $\theta=2$, that is, $S_\theta(t)=e^{t\Delta}$. 
Let $\tau=t^{1/2}$.
It follows from \eqref{eq:2.16} that 
$$
\left|\left[e^{t\Delta}\varphi\right](x)\right|
 \le C\int_{{\mathbb R}^n}t^{-\frac{n}{2}}\left(1+t^{-\frac{1}{2}}|x-y|\right)^{-n-1}|\varphi(y)|\,dy
 \le C\left[S_1(\tau)|\varphi|\right](x)
$$
for $(x,t)\in{\mathbb R}^n\times(0,\infty)$. 
This together with Proposition~\ref{Proposition:3.2} in the case of $\theta=1$ implies that 
\begin{align*}
|||e^{t\Delta}\varphi|||_{q,\beta;T^{1/2}}
 & \le C|||S_1(\tau)|\varphi||||_{q,\beta;T^{1/2}}\\
 & \le C\tau^{-n\left(\frac{1}{r}-\frac{1}{q}\right)}\left[\log\left(e+\frac{1}{\tau}\right)\right]^{-\frac{\alpha}{r}+\frac{\beta}{q}}|||\varphi|||_{r,\alpha;T^{1/2}}\\
 & \le Ct^{-\frac{n}{2}\left(\frac{1}{r}-\frac{1}{q}\right)}\left[\log\left(e+\frac{1}{t}\right)\right]^{-\frac{\alpha}{r}+\frac{\beta}{q}}|||\varphi|||_{r,\alpha;T^{1/2}}
 \quad\mbox{for}\quad t\in(0,T).
\end{align*}
Thus Proposition~\ref{Proposition:3.2} follows in the case of $\theta=2$. 
The proof of Proposition~\ref{Proposition:3.2} is complete.~$\Box$
\section{Proof of Theorems~\ref{Theorem:1.1} and \ref{Theorem:1.2}}
We apply the contraction mapping theorem to problem~\eqref{eq:P} in uniformly local {weak} Zygmund type spaces,  
and prove Theorems~\ref{Theorem:1.1} and \ref{Theorem:1.2}. 
We also prove Corollary~\ref{Corollary:1.1}. 
We first prove the following proposition.
\begin{proposition}
\label{Proposition:4.1}
Let $p=p_\theta$, $T_*\in(0,\infty)$, and $\gamma\in[0,n/\theta)$. 
Then there exists $\epsilon>0$ such that if $\varphi\in {\mathfrak L}_{{\rm ul}}^{1,\infty}(\log{\mathfrak L})^{n/\theta}$ satisfies
\begin{equation}
\label{eq:4.1}
|||\varphi|||_{1,n/\theta;T^{1/\theta}}\le\epsilon\quad\mbox{for some $T\in(0,T_*]$},
\end{equation}
then problem~\eqref{eq:P} possesses a solution 
$u\in C((0,T):{\mathfrak L}^{1,\infty}_{{\rm ul}}(\log{\mathfrak L})^{n/\theta})\cap L^\infty_{\rm loc}(0,T:L^\infty)$ 
in ${\mathbb R}^n\times(0,T)$, 
with~$u$ satisfying 
\begin{equation}
\label{eq:4.2}
	\begin{aligned}
 	& 
	|||u(t)|||_{1,n/\theta;T^{1/\theta}}\le C|||\varphi|||_{1,n/\theta;T^{1/\theta}},
	\\
 	& 
	|||u(t)|||_{p,\gamma;T^{1/\theta}}
 	\le Ct^{-\frac{n}{\theta}\left(1-\frac{1}{p}\right)}
	\left[\log\left(e+\frac{1}{t}\right)\right]^{-\frac{n}{\theta}+\frac{\gamma}{p}}|||\varphi|||_{1,n/\theta;T^{1/\theta}},
	\\
	&
	\|u(t)\|_{L^\infty}
 	\le Ct^{-\frac{n}{\theta}}
	\left[\log\left(e+\frac{1}{t}\right)\right]^{-\frac{n}{\theta}}|||\varphi|||_{1,n/\theta;T^{1/\theta}},
	\end{aligned}
\end{equation}
for $t\in(0,T)$. Here $C$ is a positive constant depending only on $T_*$, $n$, $\theta$, and $\gamma$.
\end{proposition}
Throughout this section, we set 
$$
T_*\in(0,\infty),\quad T\in(0,T_*],\quad
p:=p_\theta= 1+\frac{\theta}{n},\quad \alpha:=\frac{n}{\theta},\quad 0\le \gamma<\alpha,
\quad\varphi\in {\mathfrak L}^{1,\infty}_{{\rm ul}}(\log{\mathfrak L})^\alpha.
$$
Let $\epsilon>0$, and assume \eqref{eq:4.1}. 
By Proposition~\ref{Proposition:3.2} we find $C_*>0$ such that
\begin{equation}
\label{eq:4.3}
	\begin{aligned}
 	& 
	\sup_{t>0}\,|||S_\theta(t)\varphi|||_{1,\alpha;T^{1/\theta}}\le C_*|||\varphi|||_{1,\alpha;T^{1/\theta}}\le C_*\epsilon,
	\\
 	& 
	\sup_{t>0}\,t^{\frac{n}{\theta}\left(1-\frac{1}{p}\right)}
	\left[\log\left(e+\frac{1}{t}\right)\right]^{-\frac{\gamma}{p}+\alpha}
 	|||S_\theta(t)\varphi|||_{p,\gamma;T^{1/\theta}}\le C_*|||\varphi|||_{1,\alpha;T^{1/\theta}}\le C_*\epsilon,
	\\
	& 
	\sup_{t>0}\,t^{\frac{n}{\theta}}
	\left[\log\left(e+\frac{1}{t}\right)\right]^{\alpha}
 	\|S_\theta(t)\varphi\|_{L^\infty}\le C_*|||\varphi|||_{1,\alpha;T^{1/\theta}}\le C_*\epsilon.
	\end{aligned}
\end{equation}
Define 
$$
	X_T
	:=C((0,T):{\mathfrak L}_{{\rm ul}}^{1,\infty}(\log{\mathfrak L})^\alpha)
	\cap L^\infty_{\rm loc}((0,T):{\mathfrak L}_{{\rm ul}}^{p,\infty}(\log{\mathfrak L})^{\gamma})
	\cap L^\infty_{\rm loc}((0,T):L^\infty).
$$
Setting $C^*=2C_*$, 
for any $u\in X_T$, we say that $u\in X_T(C^*\epsilon)$ if $u$ satisfies 
\begin{equation}
\label{eq:4.4}
	\begin{aligned}
 	& 
	\sup_{0<t<T}\,|||u(t)|||_{1,\alpha;T^{1/\theta}}
	+
	\sup_{0<t<T}\,t^{\frac{n}{\theta}\left(1-\frac{1}{p}\right)}
	\left[\log\left(e+\frac{1}{t}\right)\right]^{-\frac{\gamma}{p}+\alpha}|||u(t)|||_{p,\gamma;T^{1/\theta}}
	\\
 	&\hspace{5cm}
	+
	\sup_{0<t<T}\,t^{\frac{n}{\theta}}
	\left[\log\left(e+\frac{1}{t}\right)\right]^{\alpha}\|u(t)\|_{L^\infty}
	\le C^*\epsilon.
	\end{aligned}
\end{equation}
For any $u$, $v\in X_T(C^*\epsilon)$, set
$$
	d_X(u,v) :=d^1_X(u,v)+d^2_X(u,v)+d^3_X(u,v),
$$
where
\begin{align*}
	d^1_X(u,v) 
	&
	 :=\sup_{0<t<T}\,|||u(t)-v(t)|||_{1,\alpha;T^{1/\theta}},
	 \\
	d^2_X(u,v) 
	& 
	:=\sup_{0<t<T}\,t^{\frac{n}{\theta}\left(1-\frac{1}{p}\right)}
	\left[\log\left(e+\frac{1}{t}\right)\right]^{-\frac{\gamma}{p}+\alpha}|||u(t)-v(t)|||_{p,\gamma;T^{1/\theta}},
	\\
	 d^3_X(u,v) 
	& 
	:=\sup_{0<t<T}\,t^{\frac{n}{\theta}}
	\left[\log\left(e+\frac{1}{t}\right)\right]^{\alpha}\|u(t)-v(t)\|_{L^\infty}.
\end{align*}
Then $(X_T,d_X)$ is a Banach space and $X_T(C^*\epsilon)$ is closed in $(X_T,d_X)$. 
Define 
$$
\Phi(u):=S_\theta(t)\varphi+\int_0^t S_\theta(t-s)F_p(u(s))\,ds\quad\mbox{for}\quad u\in X_T(C^*\epsilon),
$$
{where $F_p(s)=|s|^{p-1}s$ for $s\in \mathbb{R}$.}
For the proof of Proposition~\ref{Proposition:4.1} we prepare the following two lemmas.
\begin{lemma}
\label{Lemma:4.1}
Let $\epsilon>0$, and assume that \eqref{eq:4.1} holds for some $T\in(0,T_*]$. 
Then there exists $C=C(n,\theta,C_*,T_*)>0$ such that
$$
d^1_X(\Phi(u),\Phi(v))+d^2_X(\Phi(u),\Phi(v))\le C\epsilon^{p-1}\,d^2_X(u,v)
\quad\mbox{for}\quad u,v\in X_T(C^*\epsilon).
$$
\end{lemma}
{\bf Proof.}
Let $u$, $v\in X_T(C^*\epsilon)$. Let $0<s<t<T$. 
It follows that 
\begin{equation}
\label{eq:4.5}
|F_p(u(x,s))-F_p(v(x,s))|\le w(x,s)|u(x,s)-v(x,s)|
\quad\mbox{for}\quad x\in{\mathbb R}^n, 
\end{equation}
where $w(x,s):=p(|u(x,s)|^{p-1}+|v(x,s)|^{p-1})$. 
Then, by Lemmas~\ref{Lemma:2.2} and \ref{Lemma:2.3} we have 
\begin{equation}
\label{eq:4.6}
\begin{split}
 & |||F_p(u(s))-F_p(v(s))|||_{1,\gamma;T^{1/\theta}}\\
 & \le |||w(s)|||_{p/(p-1),\gamma;T^{1/\theta}} |||u(s)-v(s)|||_{p,\gamma;T^{1/\theta}}\\
 & \le p\left(|||u(s)|||^{p-1}_{p,\gamma;T^{1/\theta}}+|||v(s)|||^{p-1}_{p,\gamma;T^{1/\theta}}\right)|||u(s)-v(s)|||_{p,\gamma;T^{1/\theta}}.
\end{split}
\end{equation}
Since $u$, $v\in X_T(C^*\epsilon)$, 
by \eqref{eq:4.4} we obtain 
\begin{equation}
\label{eq:4.7}
\begin{split}
 & |||F_p(u(s))-F_p(v(s))|||_{1,\gamma;T^{1/\theta}}\\
 & \le Cs^{-\frac{n(p-1)}{\theta}\left(1-\frac{1}{p}\right)}\left[\log\left(e+\frac{1}{s}\right)\right]^{\frac{\gamma(p-1)}{p}-\alpha(p-1)}
 (C^*\epsilon)^{p-1}\\
 & \qquad\quad 
 \times C s^{-\frac{n}{\theta}\left(1-\frac{1}{p}\right)}\left[\log\left(e+\frac{1}{{s}}\right)\right]^{\frac{\gamma}{p}-\alpha} d^2_X(u,v)\\
 & =C\epsilon^{p-1} s^{-\frac{n(p-1)}{\theta}}\left[\log\left(e+\frac{1}{{s}}\right)\right]^{\gamma-\alpha p}d^2_X(u,v).
\end{split}
\end{equation}
This together with Proposition~\ref{Proposition:3.2} implies that
\begin{equation}
\label{eq:4.8}
	\begin{aligned}
 	& 
	\biggr|\biggr|\biggr|\int_0^t S_\theta(t-s)[F_p(u(s))-F_p(v(s))]\,ds\,\biggr|\biggr|\biggr|_{q,\beta;T^{1/\theta}}
	\\
 	& 
	\le\int_0^t |||S_\theta(t-s)[F_p(u(s))-F_p(v(s))]|||_{q,\beta;T^{1/\theta}}\,ds
	\\
 	& \le C\int_0^t(t-s)^{-\frac{n}{\theta}\left(1-\frac{1}{q}\right)}
	\left[\log\left(e+\frac{1}{t-s}\right)\right]^{-\gamma+\frac{\beta}{q}}
	|||F_p(u(s))-F_p(v(s))|||_{1,\gamma;T^{1/\theta}}\,ds
	\\
 	& \le C\epsilon^{p-1}d^2_X(u,v)
	\\
 	& \qquad
 	\times\int_0^t (t-s)^{-\frac{n}{\theta}\left(1-\frac{1}{q}\right)}
	\left[\log\left(e+\frac{1}{t-s}\right)\right]^{-\gamma+\frac{\beta}{q}}s^{-\frac{n}{\theta}(p-1)}
	\left[\log\left(e+\frac{1}{s}\right)\right]^{\gamma-\alpha p}\,ds
	\end{aligned}
\end{equation}
for $q\in[1,p]$ and $\beta\in[\gamma,\alpha]$.

On the other hand, since
\begin{equation}
\label{eq:4.9}
\gamma-\alpha p=\gamma-\frac{n}{\theta}\left(1+\frac{\theta}{n}\right)=\gamma-\frac{n}{\theta}-1
=\gamma-\alpha-1<-1,
\end{equation}
by Lemma~\ref{Lemma:3.1}~(2) and \eqref{eq:3.1} we have
\begin{equation}
\label{eq:4.10}
	\begin{aligned}
 	& 
	\int_0^{t/2} (t-s)^{-\frac{n}{\theta}\left(1-\frac{1}{q}\right)}
	\left[\log\left(e+\frac{1}{t-s}\right)\right]^{-\gamma+\frac{\beta}{q}}s^{-\frac{n}{\theta}(p-1)}
	\left[\log\left(e+\frac{1}{s}\right)\right]^{\gamma-\alpha p}\,ds\\
 	& 
	\le Ct^{-\frac{n}{\theta}\left(1-\frac{1}{q}\right)}
	\left[\log\left(e+\frac{1}{t}\right)\right]^{-\gamma+\frac{\beta}{q}}
	\int_0^{t/2}s^{-1}\left[\log\left(e+\frac{1}{s}\right)\right]^{\gamma-\alpha p}\,ds
	\\
 	& 
	\le Ct^{-\frac{n}{\theta}\left(1-\frac{1}{q}\right)}
	\left[\log\left(e+\frac{1}{t}\right)\right]^{-\gamma+\frac{\beta}{q}}
	\cdot C\left[\log\left(e+\frac{1}{t}\right)\right]^{\gamma-\frac{n}{\theta}}
	\\
 	& 
	=Ct^{-\frac{n}{\theta}\left(1-\frac{1}{q}\right)}\left[\log\left(e+\frac{1}{t}\right)\right]^{\frac{\beta}{q}-\alpha}
	\end{aligned}
\end{equation}
for $t\in(0,T)$. 
Similarly, since
$$
	-\frac{n}{\theta}\left(1-\frac{1}{q}\right)\ge-\frac{n(p-1)}{\theta p}=-\frac{1}{p}>-1,
$$
by Lemma~\ref{Lemma:3.1}~(1) and \eqref{eq:3.1} we obtain
\begin{equation}
\label{eq:4.11}
	\begin{aligned}
 	& 
	\int_{t/2}^t (t-s)^{-\frac{n}{\theta}\left(1-\frac{1}{q}\right)}
	\left[\log\left(e+\frac{1}{t-s}\right)\right]^{-\gamma+\frac{\beta}{q}}
	s^{-\frac{n}{\theta}(p-1)}\left[\log\left(e+\frac{1}{s}\right)\right]^{\gamma-\alpha p}\,ds
	\\
 	& 
	\le Ct^{-\frac{n}{\theta}(p-1)}\left[\log\left(e+\frac{1}{t}\right)\right]^{\gamma-\alpha p}
	\int_{t/2}^t (t-s)^{-\frac{n}{\theta}\left(1-\frac{1}{q}\right)}
	\left[\log\left(e+\frac{1}{t-s}\right)\right]^{-\gamma+\frac{\beta}{q}}\,ds
	\\
 	& 
	\le Ct^{-1}\left[\log\left(e+\frac{1}{t}\right)\right]^{\gamma-\alpha p}
	\cdot Ct^{-\frac{n}{\theta}\left(1-\frac{1}{q}\right)+1}
	\left[\log\left(e+\frac{1}{t}\right)\right]^{-\gamma+\frac{\beta}{q}}
	\\
 	& 
	=Ct^{-\frac{n}{\theta}\left(1-\frac{1}{q}\right)}\left[\log\left(e+\frac{1}{t}\right)\right]^{\frac{\beta}{q}-\alpha p}
 	\le Ct^{-\frac{n}{\theta}\left(1-\frac{1}{q}\right)}\left[\log\left(e+\frac{1}{t}\right)\right]^{\frac{\beta}{q}-\alpha}
	\end{aligned}
\end{equation}
for $t\in(0,T)$. 
Combining \eqref{eq:4.8}, \eqref{eq:4.10}, and \eqref{eq:4.11} with $(q,\beta)=(1,\alpha)$ and $(p,\gamma)$, 
we deduce that
$$
	\begin{aligned}
 	& 
	d^1_X(\Phi(u),\Phi(v))+d^2_X(\Phi(u),\Phi(v))
	\\
 	&
	=\sup_{0<t<T}\,
 	\biggr|\biggr|\biggr|\int_0^t S_\theta(t-s)[F_p(u(s))-F_p(v(s))]\,ds\,\biggr|\biggr|\biggr|_{1,\alpha;T^{1/\theta}}
	\\
	&\qquad
	+\sup_{0<t<T}\,t^{\frac{n}{\theta}\left(1-\frac{1}{p}\right)}
	\left[\log\left(e+\frac{1}{t}\right)\right]^{-\frac{\gamma}{p}+\alpha}
 	\biggr|\biggr|\biggr|\int_0^t S_\theta(t-s)[F_p(u(s))-F_p(v(s))]\,ds\,\biggr|\biggr|\biggr|_{p,\gamma;T^{1/\theta}}
	\\
 	& \le C\epsilon^{p-1} d^2_X(u,v)
\end{aligned}
$$
for $u$, $v\in X_T(C^*\epsilon)$.
Thus Lemma~\ref{Lemma:4.1} follows.
$\Box$

\begin{lemma}
\label{Lemma:4.2}
Let $\epsilon>0$, and assume that \eqref{eq:4.1} holds for some $T\in(0,T_*]$. 
Then there exists $C=C(n,\theta,C_*,T_*)>0$ such that 
$$
	d^3_X(\Phi(u),\Phi(v))\le C\epsilon^{p-1}\,\bigg(d^2_X(u,v)+d^3_X(u,v)\bigg)
	\quad\mbox{for}\quad u,v\in X_T(C^*\epsilon).
$$
\end{lemma}
{\bf Proof.}
Let $u$, $v\in X_T(C^*\epsilon)$. 
Let $0<s<t<T$. 
Similarly to \eqref{eq:4.6}, we have 
$$
	\begin{aligned}
	\|F_p(u(s))-F_p(v(s))\|_{L^\infty}
 	& 
	\le \|w(s)\|_{L^\infty} \|u(s)-v(s)\|_{L^\infty}
	\\
 	& 
	\le 
	p\left(\|u(s)\|^{p-1}_{L^\infty}+\|v(s)\|^{p-1}_{L^\infty}\right)
	\|u(s)-v(s)\|_{L^\infty}.
	\end{aligned}
$$
Since $u$, $v\in X_T(C^*\epsilon)$, 
by \eqref{eq:4.4} we obtain 
$$
	\begin{aligned}
 	& 
	\|F_p(u(s))-F_p(v(s))\|_{L^\infty}
	\\
 	& 
	\le Cs^{-\frac{n(p-1)}{\theta}}\left[\log\left(e+\frac{1}{s}\right)\right]^{-\alpha(p-1)}(C^*\epsilon)^{p-1}
	 \cdot s^{-\frac{n}{\theta}}\left[\log\left(e+\frac{1}{{s}}\right)\right]^{-\alpha} d^3_X(u,v)\\
 	& 
	=C\epsilon^{p-1} s^{-\frac{np}{\theta}}\left[\log\left(e+\frac{1}{{s}}\right)\right]^{-\alpha p}d^3_X(u,v).
	\end{aligned}
$$
This together with Proposition~\ref{Proposition:3.2} and \eqref{eq:4.7} implies that 
$$
	\begin{aligned}
 	& 
	\bigg\|\int_0^t S_\theta(t-s)[F_p(u(s))-F_p(v(s))]\,ds\,\bigg\|_{L^\infty}
	\\
 	& 
	\le\int_0^t \|S_\theta(t-s)[F_p(u(s))-F_p(v(s))]\|_{L^\infty}\,ds
	\\
 	& 
	\le C\int_0^{t/2}(t-s)^{-\frac{n}{\theta}}\left[\log\left(e+\frac{1}{t-s}\right)\right]^{-\gamma}
	|||F_p(u(s))-F_p(v(s))|||_{1,\gamma;T^{1/\theta}}\,ds
	\\
	&
	\qquad\qquad\qquad
	+C\int_{t/2}^t\|F_p(u(s))-F_p(v(s))\|_{L^\infty}\,ds
	\\
 	& 
	\le C\epsilon^{p-1}d^2_X(u,v)
	\\
 	& \qquad
 	\times\int_0^{t/2} (t-s)^{-\frac{n}{\theta}}\left[\log\left(e+\frac{1}{t-s}\right)\right]^{-\gamma}
	s^{-\frac{n}{\theta}(p-1)}\left[\log\left(e+\frac{1}{s}\right)\right]^{\gamma-\alpha p}\,ds
	\\
	& \qquad\qquad
	+C\epsilon^{p-1}d^3_X(u,v)
	\int_{t/2}^ts^{-\frac{np}{\theta}}\left[\log\left(e+\frac{1}{{s}}\right)\right]^{-\alpha p}\,ds
	\\
	&
	\le C\epsilon^{p-1}t^{-\frac{n}{\theta}}\left[\log\left(e+\frac{1}{t}\right)\right]^{-\gamma}d^2_X(u,v)
	\int_0^{t/2} s^{-1}\left[\log\left(e+\frac{1}{s}\right)\right]^{\gamma-\alpha p}\,ds
	\\
	& \qquad\qquad
	+C\epsilon^{p-1}t^{-\frac{np}{\theta}+1}\left[\log\left(e+\frac{1}{{t}}\right)\right]^{-\alpha p}d^3_X(u,v).
	\end{aligned}
$$
Since $np=n+\theta$ and $\alpha p>\alpha$, we have
\begin{equation}
\label{eq:4.13}
	\begin{aligned}
	&
	\bigg\|\int_0^t S_\theta(t-s)[F_p(u(s))-F_p(v(s))]\,ds\,\bigg\|_{L^\infty}
	\\
	&
	\le C\epsilon^{p-1}t^{-\frac{n}{\theta}}\left[\log\left(e+\frac{1}{t}\right)\right]^{-\gamma}d^2_X(u,v)
	\int_0^{t/2}s^{-1}\left[\log\left(e+\frac{1}{s}\right)\right]^{\gamma-\alpha p}\,ds
	\\
	& \qquad\qquad
	+C\epsilon^{p-1}t^{-\frac{n}{\theta}}\left[\log\left(e+\frac{1}{{t}}\right)\right]^{-\alpha}d^3_X(u,v).
	\end{aligned}
\end{equation}
Furthermore, by Lemma~\ref{Lemma:3.1}~(2) and \eqref{eq:4.9} we see that
\begin{equation}
\label{eq:4.14}
	\int_0^{t/2}s^{-1}\left[\log\left(e+\frac{1}{s}\right)\right]^{\gamma-\alpha p}\,ds
	\le C\left[\log\left(e+\frac{1}{t}\right)\right]^{\gamma-\alpha}
\end{equation}
for $t\in(0,T)$. 
Combining \eqref{eq:4.13} and \eqref{eq:4.14}, 
we deduce that 
\begin{align*}
 	& 
	d^3_X(\Phi(u),\Phi(v))
	\\
 	& 
	=\sup_{0<t<T}\,t^{\frac{n}{\theta}}\left[\log\left(e+\frac{1}{t}\right)\right]^{\alpha}
	\bigg\|\int_0^t S_\theta(t-s)[F_p(u(s))-F_p(v(s))]\,ds\,\bigg\|_{L^\infty}
	\\
 	& 
	\le C\epsilon^{p-1} \bigg(d^2_X(u,v)+d^3_X(u,v)\bigg)
\end{align*}
for $u$, $v\in X_T(C^*\epsilon)$.
Thus Lemma~\ref{Lemma:4.2} follows.
$\Box$\vspace{5pt}

\noindent
{\bf Proof of Proposition~\ref{Proposition:4.1}.}
Let $T_*>0$. 
Let $\epsilon>0$ be small enough. 
Let $\varphi\in {\mathfrak L}^{1,\infty}_{{\rm ul}}(\log{\mathfrak L})^\alpha$ be such that 
$|||\varphi|||_{1,\alpha;T^{1/\theta}}<\epsilon$ for some $T\in(0,T_*]$. 
By \eqref{eq:4.3}, \eqref{eq:4.4}, and Lemma~\ref{Lemma:4.1} we have
\begin{equation}
\label{eq:4.15}
	\begin{aligned}
	&
	\sup_{t\in(0,T)}|||\Phi(u(t))|||_{1,\alpha;T^{1/\theta}}
	+
	\sup_{0<t<T}\,\left\{t^{\frac{n}{\theta}\left(1-\frac{1}{p}\right)}
	\left[\log\left(e+\frac{1}{t}\right)\right]^{-\frac{\gamma}{p}+\alpha}|||\Phi(u(t))|||_{p,\gamma;T^{1/\theta}}\right\}
	\\
	& \le |||S_\theta(t)\varphi|||_{1,\alpha,T^{1/\theta}}
	+
	\sup_{0<t<T}\,\left\{t^{\frac{n}{\theta}\left(1-\frac{1}{p}\right)}
	\left[\log\left(e+\frac{1}{t}\right)\right]^{-\frac{\gamma}{p}+\alpha}
	|||S_\theta(t)\varphi|||_{p,\gamma;T^{1/\theta}}\right\}
	\\
	&\qquad
	+d_X^1(\Phi(u),\Phi(0))+d_X^2\bigl(\Phi(u),\Phi(0)\bigr)
	\\
 	& \le C_*\epsilon+C\epsilon^{p-1}d_X^2(u,0)
 	\le C_*\epsilon+C\epsilon^{p-1}\cdot C^*\epsilon\le C^*\epsilon
	\end{aligned}
\end{equation}
for $u\in X_T(C^*\epsilon)$. 
Similarly, we observe from Lemma~\ref{Lemma:4.2}, \eqref{eq:4.3}, and \eqref{eq:4.4} that
\begin{equation}
\label{eq:4.16}
	\begin{aligned}
 	& 
	\sup_{0<t<T}\,\left\{t^{\frac{n}{\theta}}\left[\log\left(e+\frac{1}{t}\right)\right]^{\alpha}
	\|\Phi(u(t))\|_{L^\infty}\right\}
	\\
 	& 
	\le \sup_{0<t<T}\,\left\{t^{\frac{n}{\theta}}\left[\log\left(e+\frac{1}{t}\right)\right]^{\alpha}
	\|S_\theta(t)\varphi\|_{L^\infty}\right\}
 	+d_X^3(\Phi(u),\Phi(0))
	\\
 	& 
	\le C_*\epsilon+C\epsilon^{p-1}\bigg(d_X^2(u,0)+d_X^3(u,0)\bigg)
	\le C_*\epsilon+C\epsilon^{p-1}\cdot 2C^*\epsilon\le C^*\epsilon
	\end{aligned}
\end{equation}
for $u\in X_T(C^*\epsilon)$. 
By \eqref{eq:4.15} and \eqref{eq:4.16} we see that $\Phi(u)\in X_T(C^*\epsilon)$ for $u\in X_T(C^*\epsilon)$.
Furthermore, 
taking small enough $\epsilon>0$ if necessary, 
by Lemmas~\ref{Lemma:4.1} and \ref{Lemma:4.2} we have 
\begin{align*}
	d_X(\Phi(u),\Phi(v)) 
	& 
	=d^1_X(\Phi(u),\Phi(v))+d^2_X(\Phi(u),\Phi(v))+d^3_X(\Phi(u),\Phi(v))
	\\
 	& 
	\le C\epsilon^{p-1}\,\bigg(d^2_X(u,v)+d^3_X(u,v)\bigg)\le\frac{1}{2}d_X(u,v)
\end{align*}
for $u$, $v\in X_T(C^*\epsilon)$. 
Then we apply the contraction mapping theorem to find a unique $u_*\in X_T(C^*\epsilon)$ such that 
$\Phi(u_*)=u_*$ in $X_T(C^*\epsilon)$. The function~$u_*$ is a solution to problem~\eqref{eq:P} in ${\mathbb R}^n\times(0,T)$, 
with $u_*$ satisfying \eqref{eq:4.2}. 
Thus Proposition~\ref{Proposition:4.1} follows.
$\Box$\vspace{5pt}

\noindent
{\bf Proof of Theorem~\ref{Theorem:1.1}.}
Let $T>0$. Let $\varphi\in {\mathfrak L}^{1,\infty}_{{\rm ul}}(\log{\mathfrak L})^\alpha$ be such that $|||\varphi|||_{1,\alpha;T^{1/\theta}}$ is small enough. 
Then, by Proposition~\ref{Proposition:4.1} 
we find a solution $u$ to problem~\eqref{eq:P} in ${\mathbb R}^n\times(0,T)$, with $u$ satisfying \eqref{eq:4.2}. 
Let $\beta\in(\gamma,n/\theta)$. 
Then, by Proposition~\ref{Proposition:3.2}, Lemma~\ref{Lemma:2.3}, and \eqref{eq:4.2} we obtain
\begin{equation}
\label{eq:4.17}
\begin{split}
 & |||u(t)-S_\theta(t)\varphi|||_{1,\beta;T^{1/\theta}}\\
 & \le\int_0^t |||S_\theta(t-s)F_p(u(s))|||_{1,\beta;T^{1/\theta}}\,ds\\
 & \le C\int_0^t \left[\log\left(e+\frac{1}{t-s}\right)\right]^{-\gamma+\beta}|||F_p(u(s))|||_{1,\gamma;T^{1/\theta}}\,ds\\
 & =C\int_0^t \left[\log\left(e+\frac{1}{t-s}\right)\right]^{-\gamma+\beta}|||u(s)|||^p_{p,\gamma;T^{1/\theta}}\,ds\\
 & \le C|||\varphi|||_{1,\alpha;T^{1/\theta}}^p
 \int_0^t \left[\log\left(e+\frac{1}{t-s}\right)\right]^{-\gamma+\beta}s^{-1}
 \left[\log\left(e+\frac{1}{s}\right)\right]^{\gamma-\alpha p}\,ds
\end{split}
\end{equation}
for $t\in(0,T)$. 
On the other hand, since $\beta<\theta/n$, by Lemma~\ref{Lemma:3.1}~(2) and \eqref{eq:4.9} we have
\begin{equation}
\label{eq:4.18}
\begin{split}
 & \int_0^{t/2} \left[\log\left(e+\frac{1}{t-s}\right)\right]^{-\gamma+\beta}s^{-1}
\left[\log\left(e+\frac{1}{s}\right)\right]^{\gamma-\alpha p}\,ds\\
 & \le C\left[\log\left(e+\frac{1}{t}\right)\right]^{-\gamma+\beta}
\int_0^{t/2}s^{-1}\left[\log\left(e+\frac{1}{s}\right)\right]^{\gamma-\alpha p}\,ds\\
 & \le C\left[\log\left(e+\frac{1}{t}\right)\right]^{-\gamma+\beta}\cdot C\left[\log\left(e+\frac{1}{t}\right)\right]^{\gamma-\frac{n}{\theta}}\to 0
\end{split}
\end{equation}
and 
\begin{equation}
\label{eq:4.19}
\begin{split}
 & \int_{t/2}^t \left[\log\left(e+\frac{1}{t-s}\right)\right]^{-\gamma+\beta}s^{-1}
\left[\log\left(e+\frac{1}{s}\right)\right]^{\gamma-\alpha p}\,ds\\
 & \le Ct^{-1}\left[\log\left(e+\frac{1}{t}\right)\right]^{\gamma-\frac{n}{\theta}-1}
\int_{t/2}^t\left[\log\left(e+\frac{1}{t-s}\right)\right]^{-\gamma+\beta}\,ds\\
 & \le Ct^{-1}\left[\log\left(e+\frac{1}{t}\right)\right]^{\gamma-\frac{n}{\theta}-1}\cdot 
 Ct\left[\log\left(e+\frac{1}{t}\right)\right]^{-\gamma+\beta}\to 0
\end{split}
\end{equation}
as $t\to +0$. 
Combining \eqref{eq:4.17}, \eqref{eq:4.18}, and \eqref{eq:4.19}, we see that 
$$
\lim_{t\to +0}|||u(t)-S_\theta(t)\varphi|||_{1,\beta;T^{1/\theta}}=0\quad\mbox{for $\beta\in(\gamma,n/\theta)$}. 
$$
This together with \eqref{eq:2.12} implies that 
\begin{equation}
\label{eq:4.20}
\lim_{t\to +0}|||u(t)-S_\theta(t)\varphi|||_{1,\beta;T^{1/\theta}}=0\quad\mbox{for $\beta\in[0,n/\theta)$}. 
\end{equation}

It remains to prove that $u\to \varphi$ in the sense of distributions. 
Let $\eta\in C_0({\mathbb R}^n)$. Let $R>0$ be such that $\mbox{supp}\,\eta\subset B(0,R)$. 
By {\eqref{eq:1.3}}, \eqref{eq:1.4}, and \eqref{eq:4.20} we have
\begin{equation}
\label{eq:4.21}
\begin{split}
\left|\int_{{\mathbb R}^n}
\left(u(x,t)-[S_\theta(t)\varphi](x)\right)\eta(x)\,dx\right|
 & \le C\|\eta\|_{L^\infty}\int_{B(0,R)}|u(x,t)-[S_\theta(t)\varphi](x)|\,dx\\
 & \le C\|\eta\|_{L^\infty}|||u(t)-S_\theta(t)\varphi|||_{1,0;T^{1/\theta}}\to 0
\end{split}
\end{equation}
as $t\to+0$. 
Set 
$$
\eta(x,t):=\int_{{\mathbb R}^n}G_\theta(x-y,t)\eta(y)\,dy
\quad\mbox{for}\quad (x,t)\in{\mathbb R}^n\times(0,\infty).
$$
It follows from \eqref{eq:2.18} that 
\begin{equation}
\label{eq:4.22}
\lim_{t\to +0}\|{\eta(\cdot,t)}-\eta\|_{L^\infty}=0. 
\end{equation}
On the other hand, by \eqref{eq:2.16} we have
\begin{align*}
|\eta(x,t)| & \le Ct^{-\frac{n}{\theta}}\int_{B(0,R)}(1+t^{-\frac{1}{\theta}}|x-y|)^{-n-\theta}|\eta(y)|\,dy\\
 & \le C\|\eta\|_{L^\infty}t^{-\frac{n}{\theta}}\cdot C(t^{-\frac{1}{\theta}}|x|)^{-n-\theta}
\le T|x|^{-n-\theta}
\end{align*}
for $x\in{\mathbb R}^n\setminus B(0,2R)$ and $t\in(0,T)$. 
Since $\|{\eta(\cdot,t)}\|_{L^\infty}\le \|\eta\|_{L^\infty}$ for $t>0$, 
we obtain
\begin{equation}
\label{eq:4.23}
|\eta(x,t)|\le C(1+|x|)^{-n-\theta}\quad\mbox{for}\quad (x,t)\in{\mathbb R}^n\times(0,T).
\end{equation}
Furthermore, it follows from Proposition~\ref{Proposition:3.2} with $q=\infty$ that 
$$
[S_\theta(1)|\varphi|](0)=\int_{{\mathbb R}^n}G_\theta(y,1)|\varphi(y)|\,dy<\infty. 
$$
This together with \eqref{eq:2.16} implies that 
\begin{equation}
\label{eq:4.24}
\int_{{\mathbb R}^n}(1+|y|)^{-n-\theta}|\varphi(y)|\,dy<\infty.
\end{equation}
Therefore, by \eqref{eq:4.22}, \eqref{eq:4.23}, and \eqref{eq:4.24} 
we apply the Fubini theorem and the Lebesgue convergence theorem to obtain 
\begin{align*}
\int_{{\mathbb R}^n}[S_\theta(t)\varphi](x)\eta(x)\,dx
 & =\int_{{\mathbb R}^n}\left(\int_{{\mathbb R}^n}G_\theta(x-y,t)\varphi(y)\,dy\right)\eta(x)\,dx\\
 & =\int_{{\mathbb R}^n}\left(\int_{{\mathbb R}^n}G_\theta(x-y,t)\eta(x)\,dx\right)\varphi(y)\,dy
=\int_{{\mathbb R}^n}\eta(y,t)\varphi(y)\,dy\\
 & \to \int_{{\mathbb R}^n}\eta(y)\varphi(y)\,dy
\end{align*}
as $t\to +0$. 
Then we deduce from \eqref{eq:4.21} that 
$$
\lim_{t\to +0}\int_{{\mathbb R}^n}u(x,t)\eta(x)\,dx=\int_{{\mathbb R}^n}\varphi(x)\eta(x)\,dx
\quad\mbox{for}\quad\eta\in C_0({\mathbb R}^n),
$$
that is, $u(t)\to \varphi$ in the sense of distributions. 
The proof of Theorem~\ref{Theorem:1.1} is complete.
$\Box$\vspace{5pt}
\newline
{\bf Proof of Corollary~\ref{Corollary:1.1}.} 
Let $\varphi_c$ be as in \eqref{eq:1.1} with $p=p_\theta$. 
It follows from the definition of the non-increasing rearrangements that 
\begin{equation}
\label{eq:4.25}
(\varphi_c)^*(s)\le Cs^{-1}\left[\log\left(e+\frac{1}{s}\right)\right]^{-\frac{n}{\theta}-1}\quad\mbox{for}\quad s\in(0,\infty).
\end{equation}
Let $S>0$.
Then, by Lemma~\ref{Lemma:3.1}~(2), \eqref{eq:2.3}, and \eqref{eq:4.25} we see that
$$
(\varphi_c)^{**}(s)\le Cs^{-1}\left[\log\left(e+\frac{1}{s}\right)\right]^{-\frac{n}{\theta}}\quad\mbox{for}\quad s\in(0,S).
$$
This implies that $\varphi_c\in {\mathfrak L}^{1,\infty}_{{\rm ul}}(\log{\mathfrak L})^{n/\theta}$. 
Then Corollary~\ref{Corollary:1.1} follows from Theorem~\ref{Theorem:1.1}.
$\Box$ 
\vspace{5pt}
\newline
{\bf Proof of Theorem~\ref{Theorem:1.2}.}
Since $\alpha>n/\theta$, 
it follows from \eqref{eq:2.12} that 
$$
|||\varphi|||_{1,n/\theta;T^{1/\theta}}
\le C\left[\log\left(e+\frac{1}{T^{1/\theta}}\right)\right]^{\frac{n}{\theta}-\alpha}|||\varphi|||_{1,\alpha;T^{1/\theta}}\to 0
\quad\mbox{as}\quad T\to +0.
$$
Then, by Theorem~\ref{Theorem:1.1} we find a solution $u$ to problem~\eqref{eq:P} in ${\mathbb R}^n\times(0,T)$ 
for some small enough $T>0$, with $u$ satisfying \eqref{eq:1.6} and \eqref{eq:1.7}. 
Thus Theorem~\ref{Theorem:1.2} follows.
$\Box$
\vspace{5pt}

At the end of this paper 
we recall the definitions of the usual Zygmund space and the usual weak Zygmund space, 
and explain the advantage of our weak Zygmund type spaces.  
\begin{remark}
\label{Remark:4.1}
{\rm (i)} 
We recall the Zygmund space $L^q(\log L)^\alpha$ and the weak Zygmund space $L^{q,\infty}(\log L)^\alpha$. 
For any $q\in[1,\infty]$ and $\alpha\ge 0$, set 
\begin{align*}
L^q(\log L)^\alpha & :=\{f\in L^1_{\rm loc}(\mathbb R^n)\,:\, \|f\|_{L^q(\log L)^\alpha}<\infty\},\\
L^{q,\infty}(\log L)^\alpha & :=\{f\in L^1_{\rm loc}(\mathbb R^n)\,:\, \|f\|_{L^{q,\infty}(\log L)^\alpha}<\infty\},
\end{align*}
where
\begin{align}
\label{eq:4.26}
	\|f\|_{L^q(\log L)^\alpha}
	 & :=\left(\int_0^\infty \bigg[\log\left(e+\frac{1}{s}\right)\bigg]^\alpha
	f^*(s)^q\,ds\right)^\frac{1}{q},\\
	\label{eq:4.27}
	\|f\|_{L^{q,\infty}(\log L)^\alpha}
	 & :=\sup_{s>0}\,
	\left\{\left[\log\left(e+\frac{1}{s}\right)\right]^\alpha sf^*(s)^q\right\}^{\frac{1}{q}}.
\end{align}
See e.g., \cite{BS}*{\it{Chapter~$4$, Section~$6$}} and \cite{Wadade}.
For the case $q>1$,
as in the Lorentz space {\rm({\it see e.g., \cite{Grafakos}*{\it{Chapter~$1$, Exercises~$1.4.3$}}})},
applying Hardy's inequality {\rm ({\it see Lemma}~\ref{Lemma:2.4})} and Lemma~{\rm\ref{Lemma:3.1}} with \eqref{eq:2.3},
for any $f\in L^1_{{\rm loc}}$,  we see that $f\in L^q(\log L)^\alpha$ if and only if 
$$
	[f]_{L^q(\log L)^\alpha}
	:=\left(\int_0^\infty \bigg[\log\left(e+\frac{1}{s}\right)\bigg]^\alpha
	f^{**}(s)^q\,ds\right)^\frac{1}{q}<\infty.
$$
In contrast,
the above relation does not hold for the case $q=1$.
In fact, applying the integration by parts, we see that
$$
	\int_0^\infty \bigg[\log\left(e+\frac{1}{s}\right)\bigg]^\alpha f^*(s)\,ds
	=\alpha\int_0^\infty \bigg[\log\left(e+\frac{1}{s}\right)\bigg]^{\alpha-1}f^{**}(s)\frac{ds}{es+1}
	+\|f\|_{L^1}.
$$
{\rm (ii)} 
By O'Neil's inequality~\eqref{eq:2.7} 
we have the inequality 
$$
\left(G_\theta(\cdot,t)*\varphi\right)^{**}(s)\le\int_s^\infty \left(G_\theta(\cdot,t)\right)^{**}(\tau)\varphi^{**}(\tau)\,d\tau,
\quad s>0,
$$
which is crucial in the proof of our sharp decay estimates of $S_\theta(t)\varphi$. 
Our Zygmund type spaces are defined by the average of the non-increasing rearrangement,  
and they are effectively used in the proof of our sharp decay estimates of $S_\theta(t)\varphi$ {\rm({\it see the proof of Proposition}~\ref{Proposition:3.1})}. 
These sharp decay estimates of $S_\theta(t)\varphi$ in the spaces $\mathfrak L^{q,\infty}(\log \mathfrak L)^\alpha$
enable us to obtain Theorem~{\rm\ref{Theorem:1.1}}.

On the other hand, 
since the weak Zygmund space $L^{q,\infty}(\log L)^\alpha$ is defined by the non-increasing rearrangement,  
the inequality 
\begin{equation}
\label{eq:4.28}
\left(G_\theta(\cdot,t)*\varphi\right)^{*}(s)\le
\left(G_\theta(\cdot,t)*\varphi\right)^{**}(s)\le\int_s^\infty \left(G_\theta(\cdot,t)\right)^{**}(\tau)\varphi^{**}(\tau)\,d\tau,
\quad s>0,
\end{equation} 
seems useful for the study of decay estimates of $S_\theta(t)\varphi$ in the space $L^{q,\infty}(\log L)^\alpha$.
The first inequality in \eqref{eq:4.28} follows from inequality \eqref{eq:2.4}. 
However, in general, inequality \eqref{eq:2.4} is not sharp in $L^{1,\infty}(\log L)^\alpha$, where $\alpha>1$. 
Indeed, let $f\in L^1_{{\rm loc}}$ be such that 
$$
f^*(s)=s^{-1}\left[\log\left(e+\frac{1}{s}\right)\right]^{-\alpha},\quad s>0,
$$
where $\alpha>1$. 
Then $f\in L^{1,\infty}(\log L)^\alpha$ and 
$$
f^{**}(s)\asymp s^{-1}\left[\log\left(e+\frac{1}{s}\right)\right]^{-\alpha+1}
$$
for small enough $s>0$. 
Then $f^*(s)/f^{**}(s)\to 0$ as $s\to+0$, and we see that inequality \eqref{eq:2.4} is not sharp. 
This suggests that it is difficult to obtain sharp decay estimates of $S_\theta(t)\varphi$ in the usual weak Zygmund spaces.
\vspace{3pt}
\newline
{\rm (iii)} 
In order to overcome the disadvantage of the usual weak Zygmund spaces, 
one might consider the following weak Zygmund type spaces
$$
{\mathbb L}^{q,\infty}(\log {\mathbb L})^\alpha:=\left\{f\in L^1_{{\rm loc}}\,:\,\|f\|_{{\mathbb L}^{q,\infty}(\log {\mathbb L})^\alpha}<\infty\right\},
$$
where $1\le q<\infty$, $\alpha\ge 0$, and
$$
\|f\|_{{\mathbb L}^{q,\infty}(\log {\mathbb L})^\alpha}:=\sup_{s>0}\left\{\left[\log\left(e+\frac{1}{s}\right)\right]^{\alpha}sf^{**}(s)^q\right\}^{\frac{1}{q}}.
$$
Indeed, applying the arguments to those in the proof of Proposition~{\rm\ref{Proposition:3.2}},
we can obtain similar sharp decay estimates of $S_\theta(t)\varphi$ in the weal Zygmund type space ${\mathbb L}^{q,\infty}(\log {\mathbb L})^\alpha$ 
to those in Proposition~{\rm\ref{Proposition:3.2}}. 

On the other hand, in the proof of Theorem~{\rm\ref{Theorem:1.1}}, 
we used the inequality
\begin{equation}
\label{eq:4.29}
\| |f|^p\|_{{\mathfrak L}^{1,\infty}(\log{\mathfrak L})^\alpha}\le C\|f\|^p_{{\mathfrak L}^{p,\infty}(\log{\mathfrak L})^{\alpha}}
\quad\mbox{for}\quad f\in {\mathfrak L}^{p,\infty}(\log{\mathfrak L})^{\alpha}
\end{equation}
in order to estimate the nonlinear term $|u|^{p-1}u$, where $p>1$ and $\alpha\ge 0$. 
Actually, \eqref{eq:4.29} holds with $C=1$ and $``\le"$ replace by $``="$ {\rm ({\it see Lemma}~{\rm \ref{Lemma:2.3}})}. 
In the case of ${\mathbb L}^{q,\infty}(\log {\mathbb L})^\alpha$, 
it follows from \eqref{eq:2.5} that 
\begin{align*}
\| |f|^p\|_{{\mathbb L}^{1,\infty}(\log {\mathbb L})^\alpha}
 & =\sup_{s>0}\left\{\left[\log\left(e+\frac{1}{s}\right)\right]^{\alpha}s(|f|^p)^{**}(s)\right\}\\
 & \ge\sup_{s>0}\left\{\left[\log\left(e+\frac{1}{s}\right)\right]^{\alpha}s(f^{**}(s))^p\right\}
= \|f\|^p_{{\mathbb L}^{p,\infty}(\log {\mathbb L})^{\alpha}}
\end{align*}
for $f\in {\mathbb L}^{p,\infty}(\log {\mathbb L})^{\alpha}$, 
that is, the reverse to the desired inequality holds. 
This suggests that it is difficult to obtain a similar result to that of Theorem~{\rm\ref{Theorem:1.1}}
in the framework of weak Zygmund type spaces ${\mathbb L}^{q,\infty}(\log {\mathbb L})^\alpha$. 
\end{remark}
\appendix\section{{Appendix}}
\label{section:A}
In the Appendix we prove two propositions on relations 
between $L^q(\log L)^\alpha$, $L^{q,\infty}(\log L)^\alpha$, and $\mathfrak L^{q,\infty}(\log \mathfrak L)^\alpha$. 
We remark that 
the following relations hold for $\alpha=0$.
$$
	L^q=L^q(\log L)^0=\mathfrak L^{q,\infty}(\log \mathfrak L)^0\subsetneq L^{q,\infty}=L^{q,\infty}(\log L)^0\quad\mbox{if}\quad q\in[1,\infty).
$$
\begin{proposition}
\label{Proposition:A.1}
Let $1\le q <\infty$ and $\alpha\ge0$. 
Then
\begin{align*}
\|f\|_{\mathfrak L^{q,\infty}(\log \mathfrak L)^\alpha}
& \le 
\|f\|_{L^q(\log L)^\alpha} \quad\mbox{for}\quad f\in L^q(\log L)^\alpha,\\
 \|f\|_{L^{q,\infty}(\log L)^\alpha}
& 
\le \|f\|_{\mathfrak L^{q,\infty}(\log \mathfrak L)^\alpha}\quad\mbox{for}\quad f\in \mathfrak L^{q,\infty}(\log \mathfrak L)^\alpha.
\end{align*}
Furthermore, 
$$
	L^q(\log L)^\alpha\subsetneq \mathfrak L^{q,\infty}(\log \mathfrak L)^\alpha\subsetneq L^{q,\infty}(\log L)^\alpha,\qquad \alpha>0.
$$
\end{proposition}
\noindent
{\bf Proof.}
By \eqref{eq:2.9}, \eqref{eq:2.11}, and \eqref{eq:4.26} we see that
$$
	\begin{aligned}
	\|f\|_{\mathfrak L^{q,\infty}(\log \mathfrak L)^\alpha}
	&
	=\sup_{s>0}\,\left\{\left[\log\left(e+\frac{1}{s}\right)\right]^\alpha
	\int_0^s (f^*(\tau))^q\,d\tau\right\}^{\frac{1}{q}}
	\\
	&
	\le
	\sup_{s>0}
	\bigg(\int_0^s\left[\log\left(e+\frac{1}{\tau}\right)\right]^\alpha
	(f^*(\tau))^q\,d\tau\bigg)^\frac{1}{q}
	=\|f\|_{L^q(\log L)^\alpha}
	\end{aligned}
$$
for $f\in L^q(\log L)^{\alpha}$.
This implies that
$L^q(\log L)^\alpha\subset\mathfrak L^{q,\infty}(\log \mathfrak L)^\alpha$. 
Let $g$ be a function in ${\mathbb R}^n$ such that 
\[
g^*(s)
=\frac{d}{ds}\left\{\left[\log\left(e+\frac{1}{s}\right)\right]^{-\alpha}\right\}\chi_{(0,\delta)}(s)
=\frac{\alpha}{es^2+s}\left[\log\left(e+\frac{1}{s}\right)\right]^{-\alpha-1}\chi_{(0,\delta)}(s),
\]
where $\delta>0$ is taken so that $g^*$ is decreasing. 
Set $f(x):=|g(x)|^{\frac{1}{q}}$. 
It follows from \eqref{eq:2.1} that $f^*(s)^q=g^*(s)$.
Furthermore, 
\[
\|f\|^q_{\mathfrak L^{q,\infty}(\log \mathfrak L)^\alpha}=\sup_{s>0}\left[\log\left(e+\frac{1}{s}\right)\right]^{\alpha}\int_0^s g^*(\eta)d\eta=1
\]
and
\[
\begin{aligned}
\|f\|^q_{L^q(\log L)^{\alpha}}
=
\int_0^{\infty}\left[\log\left(e+\frac{1}{\eta}\right)\right]^{\alpha}g^*(\eta)d\eta
=
\int_0^{\delta}\frac{\alpha}{e\eta^2+\eta}\left[\log\left(e+\frac{1}{\eta}\right)\right]^{-1}d\eta
=\infty.
\end{aligned}
\]
Thus $L^q(\log L)^\alpha\subsetneq \mathfrak L^{q,\infty}(\log \mathfrak L)^\alpha.$

On the other hand, 
it follows from \eqref{eq:2.1}, \eqref{eq:2.4}, \eqref{eq:2.9}, and \eqref{eq:4.27} that
$$
	\begin{aligned}
	\|f\|_{\mathfrak L^{q,\infty}(\log \mathfrak L)^\alpha}
	&
	=\sup_{s>0}\,\bigg\{
	\left[\log\left(e+\frac{1}{s}\right)\right]^\alpha s\,(|f|^q)^{**}(s)\bigg\}^\frac{1}{q}
	\\
	&
	\ge
	\sup_{s>0}\,
	\bigg\{\left[\log\left(e+\frac{1}{s}\right)\right]^\alpha s\,(|f|^q)^*(s)\bigg\}^\frac{1}{q}
	\\
	&
	=\sup_{s>0}\,
	\bigg\{\left[\log\left(e+\frac{1}{s}\right)\right]^\alpha s\,(f^*(s))^q\bigg\}^\frac{1}{q}
	=\|f\|_{L^{q,\infty}(\log L)^\alpha}
	\end{aligned}
$$
for $f\in \mathfrak L^{q,\infty}(\log \mathfrak L)^\alpha$, and hence 
$\mathfrak L^{q,\infty}(\log \mathfrak L)^\alpha\subset L^{q,\infty}(\log L)^\alpha$.
We finally show that the inclusion is strict. 
Let $f$ be a function such that
$$
	f^*(s)=s^{-\frac{1}{q}}\left[\log\left(e+\frac{1}{s}\right)\right]^{-\frac{\alpha}{q}}\ \chi_{(0,\delta)}(s),
$$
where $\delta>0$ is chosen so that $f^*$ is decreasing.
Then $\|f\|_{L^{q,\infty}(\log L)^{\alpha}}=1$.
On the other hand,
for the case $\alpha\le 1$, we see that
$$
	s(|f|^q)^{**}(s)
	=\int_0^s \eta^{-1}\left[\log\left(e+\frac{1}{\eta}\right)\right]^{-\alpha}\,d\eta
	\ge
	\int_0^s \eta^{-1}\left[\log\left(e+\frac{1}{\eta}\right)\right]^{-1}\,d\eta=\infty
$$
for $s\in(0,\delta)$.
This implies that $f\notin \mathfrak L^{q,\infty}(\log \mathfrak L)^\alpha$.
Furthermore, for the case $\alpha>1$,
there exists $C>0$ such that 
$$
	\begin{aligned}
	s(|f|^q)^{**}(s)
	&
	=\int_0^s \eta^{-1}\left[\log\left(e+\frac{1}{\eta}\right)\right]^{-\alpha}\,d\eta
	\\
	&
	\ge 
	C\int_0^s (e\eta^2+\eta)^{-1}\left[\log\left(e+\frac{1}{\eta}\right)\right]^{-\alpha}\,d\eta
	=\frac{C}{\alpha -1}{\left[\log\left(e+\frac{1}{s}\right)\right]^{1-\alpha}}
	\end{aligned}
$$
for $s\in(0,\delta)$.
In conclusion, there exists $C>0$ such that
$$
	\|f\|_{\mathfrak L^{q,\infty}(\log \mathfrak L)^\alpha}
	\ge 
	C \sup_{0<s<\delta}\,\left[\log\left(e+\frac{1}{s}\right)\right]^{\frac{1}{q}}= \infty.
$$
Thus $\mathfrak L^{q,\infty}(\log \mathfrak L)^\alpha\subsetneq L^{q,\infty}(\log L)^\alpha$. 
The proof of Proposition~\ref{Proposition:A.1} is complete.
$\Box$\vspace{5pt}

Let $f$ be a locally integrable function in ${\mathbb R}^n$ such that 
\begin{equation}
\label{eq:A.1}
	\left[\log\left(e+\frac{1}{s}\right)\right]^{\alpha}s(|f|^q)^{**}(s)=1,\quad s>0,
\end{equation}
which is a typical function in $\mathfrak L^{q,\infty}(\log \mathfrak L)^\alpha$.
By \eqref{eq:A.1} we see 
$$
	f^*(s)^q
	=\frac{d}{ds}\left(s(|f|^q)^{**}(s)\right)
	=\frac{\alpha}{es^2+s}\left[\log\left(e+\frac{1}{s}\right)\right]^{-\alpha-1},\quad s>0. 
$$
Since
$$
	f^*(s)
	\asymp
 	s^{-\frac{1}{q}}\left[\log\left(e+\frac{1}{s}\right)\right]^{-\frac{\alpha+1}{q}}
	\quad\mbox{for small enough $s>0$},
$$
we see that $f$ also has a typical singularity of functions in $L^{q,\infty}(\log L)^{\alpha+1}$.  
These arguments suggest that $\mathfrak L^{q,\infty}(\log \mathfrak L)^\alpha$
is closely related to $L^{q,\infty}(\log L)^{\alpha+1}$.
\begin{proposition}
\label{Proposition:A.2}
Let $1\le q <\infty$ and $\alpha>0$. Then there exists $C>0$ such that 
\begin{equation}
\label{eq:A.2}
\|f\|_{\mathfrak L^{q,\infty}(\log \mathfrak L)^\alpha}\le C\|f\|_{L^{q,\infty}(\log L)^{\alpha+1}}
\end{equation}
for $f\in L^{q,\infty}(\log L)^{\alpha+1}$. 
Furthermore, 
\begin{equation}
\label{eq:A.3}
\inf\left\{\frac{\|f\|_{\mathfrak L^{q,\infty}(\log \mathfrak L)^\alpha}}{\|f\|_{L^{q,\infty}(\log L)^{\alpha+1}}}\,:\,f\in L^{q,\infty}(\log L)^{\alpha+1}\right\}=0. 
\end{equation}
\end{proposition}
{\bf Proof.}
By Lemma~\ref{Lemma:2.3} and \eqref{eq:2.1} 
it suffices to consider the case $q=1$. 
We first prove \eqref{eq:A.2} with $q=1$. 
Let $f\in L^{1,\infty}(\log L)^{\alpha+1}$, where $\alpha>0$. 
By Lemma~\ref{Lemma:3.1}~(2),
for any $R>0$, we have
$$
\begin{aligned}
	sf^{**}(s)=\int_0^s f^*(\eta)d\eta
	&
	\le
	\left(\int_0^s\eta^{-1} \left[\log\left(e+\frac{1}{\eta}\right)\right]^{-\alpha-1}d\eta
	\right)\,
	\left(
	\sup_{\eta>0}\left[\log\left(e+\frac{1}{\eta}\right)\right]^{\alpha+1}\eta f^*(\eta)
	\right)
	\\
	&
	\le
	C\left[\log\left(e+\frac{1}{s}\right)\right]^{-\alpha}
	\left(
	\sup_{\eta>0} \left[\log\left(e+\frac{1}{\eta}\right)\right]^{\alpha+1}\eta f^*(\eta)
	\right)
\end{aligned}
$$
for $s\in(0,R)$. 
This together with \eqref{eq:2.9} implies that
\begin{align*}
\|f\|_{\mathfrak L^{1,\infty}(\log \mathfrak L)^\alpha} & =
	\sup_{0<s<R}\left[\log\left(e+\frac{1}{s}\right)\right]^{\alpha}sf^{**}(s)\\
	 & \le
	C \sup_{\eta>0}\left[\log\left(e+\frac{1}{\eta}\right)\right]^{\alpha+1}\eta f^*(\eta)
	=C\|f\|_{L^{1,\infty}(\log L)^{\alpha+1}}.
\end{align*}
Thus \eqref{eq:A.2} holds for $q=1$, and the proof of \eqref{eq:A.2} is complete. 

Next, we prove \eqref{eq:A.3} with $q=1$. 
Let $\{f_n\}$ be a sequence in $L^1_{{\rm loc}}$ such that
\[
f_n^*(s)=n\left[\log(e+n)\right]^{-\alpha-1}\chi_{(0,\frac{1}{n})}(s).
\]
Since
\begin{equation}
\label{eq:A.4}
f_n^{**}(s)
=
\left\{
\begin{array}{ll}
n\left[\log(e+n)\right]^{-\alpha-1} & \mbox{for}\quad s\in\displaystyle{\left(0,\frac{1}{n}\right)},\vspace{7pt}\\
s^{-1}\left[\log(e+n)\right]^{-\alpha-1} & \mbox{for}\quad s\in\displaystyle{\left[\frac{1}{n},\infty\right)},
\end{array}
\right.
\end{equation}
we have
\begin{align*}
\|f_n\|_{\mathfrak L^{1,\infty}(\log \mathfrak L)^\alpha}
 & =\sup_{s>0}\left[\log\left(e+\frac{1}{s}\right)\right]^{\alpha}sf_n^{**}(s)\\
 & =
\left.
\left[\log\left(e+\frac{1}{s}\right)\right]^{\alpha}sf_n^{**}(s)
\right|_{s=\frac{1}{n}}=
[\log(e+n)]^{-1}
\end{align*}
for $n=1,2,\dots$. 
On the other hand, similarly to \eqref{eq:3.2}, 
we find $L\in[e,\infty)$ such that 
$$
\mbox{the function $\displaystyle{(0,\infty)\ni \tau\mapsto \tau\left[\log\left(L+\frac{1}{\tau}\right)\right]^{\alpha+1}}$ is non-decreasing}.
$$
Then, by \eqref{eq:3.1} and \eqref{eq:A.4}
we have 
\begin{align*}
 & \|f_n\|_{L^{1,\infty}(\log L)^{\alpha+1}}=\sup_{s>0}
\left[\log\left(e+\frac{1}{s}\right)\right]^{\alpha+1}sf_n^{*}(s)\\
 & \qquad
 \ge C\sup_{s>0}
\left[\log\left(L+\frac{1}{s}\right)\right]^{\alpha+1}sf_n^{*}(s)
=C\left[\log\left(L+\frac{1}{s}\right)\right]^{\alpha+1}sf_n^{*}(s)\biggr|_{s=1/n}\ge C
\end{align*}
for $n=1,2,\dots$. 
These imply \eqref{eq:A.3}. Thus Proposition~\ref{Proposition:A.2} follows. 
$\Box$
\medskip

\noindent
{\bf Acknowledgements.} 
K. I. and T. K. were supported in part by JSPS KAKENHI Grant Number JP 19H05599. 
T. K. was also supported in part by JSPS KAKENHI Grant Number JP 20K03689 and JP 22KK0035.
{N. I. was supported in part by JSPS KAKENHI Grant Number JP 21K18582.}

\end{document}